\documentclass[a4paper,11pt]{amsart}

\usepackage[T1]{fontenc}
\usepackage{lmodern}

\usepackage{xcolor}

\usepackage{amssymb,amsmath,latexsym,amsthm}
\newcommand{\RP}{\mathbb{RP}}

 \usepackage{epsfig}
  \usepackage{float}

\usepackage{graphicx}
\graphicspath{ {./illustrations/} }

\usepackage[utf8]{inputenc}
 \usepackage{url} 
 \usepackage{hyperref} 

 \usepackage{color}
 \usepackage{xcolor}
 \definecolor{tocblue}{HTML}{24476B}
 \hypersetup{	
colorlinks=true, 
linkcolor= tocblue,	
breaklinks=true, 
citecolor=blue,
urlcolor= blue, 
}

 \usepackage{url} 
 \usepackage{hyperref}

\def\C{\mathbb{C}}
\def\R{\mathbb{R}}

\def\R{\mathbb{R}}
\def\H{\mathbb{H}}
\def\Z{\mathbb{Z}}
\def\S{\mathbb{S}}

\def\F{\mathcal {F}}

\def\R{\mathbb{R}}

\def\T{\mathbb{T}}

\def\Isom{\operatorname{Isom}}

\def\Aff{\operatorname{Aff}}

\def\Id{\operatorname{Id}}
\def\PO{\operatorname{PO}}
\def\P{\operatorname{P}}

\def\Jac{\operatorname{Jac}\;}

\def\B{\mathbb{B}}

\def\\R{\mathcal R}

\def\PGL{{\operatorname{PGL}}}
\def\PSL{{\operatorname{PSL}}}

\def\O{{\operatorname{O}}}
\def\SO{{\operatorname{SO}}}
 
\def\SL{{\operatorname{SL}}}

\def\Ein{{\operatorname{Ein}} }

\def\dS{{\operatorname{dS}} }
\def\AdS{{\operatorname{AdS}} }

\def\tr{{\operatorname{tr}}}

\def\det{{\operatorname{det}}}

\def\g{{\mathfrak{g}}}

 \newcommand{\ip}[2]{\langle #1,#2\rangle}

\newtheorem{theorem}{{Theorem}}[section]
\newtheorem{proposition}[theorem]{{Proposition}}
\newtheorem{isom.ext}[theorem]{{Trivial isometric extension}}
\newtheorem{definition}[theorem]{{Definition}}
\newtheorem{lemma}[theorem]{{Lemma}}
\newtheorem{corollary}[theorem]{{Corollary}}
\newtheorem{fact}[theorem]{{\sc Fact}}
\newtheorem{remark}[theorem]{{Remark}}
\newtheorem{question}[theorem]{{Question}}
\newtheorem{example}[theorem]{{Example}}
\newtheorem{exdef}[theorem]{{Example-Definition}}

\title{Riemannian geodesic motions and null foliations }

\date{\today}

\author{S. Allout}
\address{Souheib Allout \newline
Faculty of Mathematics, USTHB University,
Algiers, Algeria}
\email{souheib.allout.math@gmail.com}

\author{A. Boukholkhal}
\address{Alaa Mohamed Boukholkhal \newline
UMPA, ENS de Lyon, France}
\email{mohamed.boukholkhal@ens-lyon.fr}

\author{L. Messaci} 

\address{Lamine Messaci \newline
Laboratoire de recherche AGM, CY Cergy-Paris Université, France }
\email{mmessaci@cyu.fr}

\author{A. Zeghib}
\address{Abdelghani Zeghib \newline
UMPA, ENS de Lyon, France }
\email{abdelghani.zeghib@ens-lyon.fr}

\begin{document}
\begin{center}
    \begin{abstract}
We study the non-collision problem: given a unit vector field $V$ on a complete Riemannian manifold $(M,g)$, under what conditions on $V$ are the maps $F_t(p)=\exp_p(tV(p))$ injective for all $t\in\R$? If such a vector field exists, we call it a geodesic motion on $M$. This occurs, for instance, when $V$ defines a geodesic foliation of $M$, in which case, the geodesic motion is called stationary. We show that geodesic motions on the Euclidean plane are necessarily stationary, that generic geodesic motions on round spheres are non-stationary, and that compact hyperbolic $3$-manifolds admit no geodesic motions. The non-collision problem has a Lorentzian interpretation: geodesic motions on $M$ correspond to one-dimensional null geodesic foliations of the spacetime $(M\times\R,g-dt^2)$. In the second part of the paper, we focus on codimension-one lightlike geodesic foliations of anti-de Sitter space $\AdS^{1,n}$. We give a complete characterization of such foliations by means of causal curves in the Einstein universe and prove rigidity in the projective anti-de Sitter case. The last part of the paper is dedicated to the study of one-dimensional geodesic foliations of $\AdS^{1,2}$. We give a complete characterization when the foliation contains at least one timelike leaf.
    \end{abstract}
\end{center}
\maketitle
\tableofcontents 

\bigskip

\tableofcontents

\section{Introduction}
Given a complete Riemannian manifold $(M,g)$ and a unit vector field $V$ on $M$, we define the following motion on $M$:
each point $p\in M$ moves along the geodesic trajectory $F_t(p)$ starting at $p$ in the direction $V(p)$. Two different moving points $p,q\in M$ meet or collide only if their trajectories intersect at the same time. Thus, the relevant notion of intersection for this motion is not simply the intersection of the corresponding geodesics in $M$ but rather the intersection of their spacetime trajectories $t\mapsto (F_t(p),t)$ in the ``spacetime'' $M\times \R$. That is, the moving points $p$ and $q$ meet if and only if there exists $t\in \R$ such that $F_t(p)=F_t(q)$.

\subsubsection{Definition}
When the motion defined by $V$ has no collisions, i.e., $F_t$ is injective for all $t$, we say that $V$ is \textbf{collision-free}.

\medskip

If $F_t$ is a homeomorphism for all $t$, the isotopy $t\mapsto F_t$ is called a {\bf geodesic motion} generated by the collision-free vector field $V$. 

Observe that, whenever globally defined, the map $F_t$ is proper. In particular, when $(M,g)$ is complete, $F_t$ is a globally defined closed map and $V$ being collision-free is equivalent to $t\mapsto F_t$ being a geodesic motion. One can ask the following natural question

\begin{question}\label{noncollisionquest}
    Characterize all geodesic motions on $(M,g)$.
\end{question}

A geodesic motion is called {\bf stationary} if its generating vector field $V$ is geodesic, $\nabla_V V = 0$ (the terms "autonomous" and "deterministic" are equally appropriate, as the discussion below explains).

\subsubsection{A spacetime interpretation: from geodesic motions to null foliations}

Asking whether a unit vector field is collision-free can be interpreted as a foliation problem. Indeed, once we add the ``time'' dimension, the trajectories $(F_t(p),t)$ are lightlike geodesics in the spacetime $(M\times\R, g-dt^2)$. The motion is collision-free whenever these trajectories define a lightlike geodesic foliation of $(M\times\R)$. Question \ref{noncollisionquest} therefore asks us to classify all lightlike geodesic foliations of $(M\times\R, g-dt^2)$.

Starting from this Riemannian-Lorentzian interplay, one is naturally led to study null one-dimensional geodesic foliations on general Lorentz manifolds, in particular the constant curvature ones. This is the main topic of the present paper.

\subsubsection{Non-autonomous equations}
A non-autonomous differential equation has the form $\dot{x} = V(t,x)$, where, for each $t$ in an interval, $x\mapsto V(t,x)$ is a vector field. A geodesic motion corresponds to the case where all solutions $t\mapsto x(t)$ of $\dot{x}(t) = V(t,x(t))$ are arc-length-parametrized geodesics. This is equivalent to $V(t,x)$ being a unit vector and $\partial_t V + \nabla_V V = 0$.
Observe, however, that this is a situation where $V(t,x)$ is fully determined by specifying $V$, i.e., $V(0,x)$, which is generally not the case for non-autonomous systems. Observe also that the spacetime construction is the same as the usual trick for making a non-autonomous system autonomous.

\subsubsection{Euler equations of fluid dynamics}
The incompressible Euler equations are classically written: 
$ \partial_t u+\nabla_u u=-\mathsf{grad}\, p,\; \mathsf{div}\,u=0$, where $u$ is a vectorial function $u: (t, x) \in U \subset \R \times \R^n \to u(t, x) \in  \R^n$, $U$ is open in $\R \times \R^n$.
In the case of constant pressure $p$, this gives rise to a volume-preserving geodesic motion (up to assuming $u$ unitary, and completeness issues), which is thus a constant-pressure solution of the Euler equations. These equations obviously can be defined for Riemannian manifolds.

\medskip

Even though this is not the objective of the present paper, let us mention the following result illustrating the interplay between geodesic motions and solutions of the Euler equations. Its proof is the same as that of the particular case given in Proposition \ref{Volume_preservation}:

\begin{theorem}\label{symmetric_divergence_free}
A geodesic motion of a compact locally symmetric space of non-positive curvature is volume preserving.
\end{theorem}

\subsection{Euclidean plane vs Minkowski space of dimension 2+1}

When $M$ is the Euclidean plane $\R^2$, the spacetime to consider is the Minkowski space $\R^{1,2}$, and the trajectories of moving points are lightlike affine lines.  

\begin{theorem}\label{thm_euclidean_motion_intro}

Let $V:\R^2\to\mathbb{S}^1$ be a $C^0$ unit vector field. Then $V$ is collision-free if and only if it constant. In particular,  $C^0$-lightlike foliations of $\R^{1,2}$ are rigid.

\end{theorem}

\begin{remark} Foliations of $\R^3$ by affine lines have been extensively studied in \cite{Salvai2009, Harrison2016, Harrison2021}. Lightlike foliations are special line foliations of $\R^3$ with allowed directions defined by the intersection of the isotropic cone with the sphere $\S^2$. The theorem above suggests that one perhaps can obtain rigidity in similar cases for foliations with restricted directions (a curve in $\S^2$ for instance).

\end{remark}

\subsubsection{Geodesic motions on flat 3-manifolds}
The preceding global rigidity on $\R^2$ is obviously equivalent to the statement that every geodesic motion is stationary. This is no longer true in dimension 3. The following is a geodesic motion on $\R^3$ which descends to the torus $\T^3 = \R^3/\Z^3$:
\[
 V(x,y,z)=\frac1{\sqrt2}
 \left(\cos z, \sin z, 1\right).
\]

Observe that $V$ is divergence-free, as stated in Theorem \ref{symmetric_divergence_free}.

\subsection{The sphere vs Einstein universe}
Hopf fibrations are one-dimensional geodesic foliations, that is, stationary geodesic motions, on odd-dimensional spheres $\S^n$. There are many works, old and recent, on this subject \cite{gluck1983great, GluckWarnerYang1983, GluckWarnerZiller1986, CahnGluckNuchi2018}. One then asks whether there are non-stationary geodesic motions. They correspond to null foliations of the Einstein universe $\Ein^{1+n} = (\S^1,-dt^2)\times\S^n$.

\begin{theorem}\label{SphereGeoMotion}
A unit vector field $V$ on $\S^n$ defines a geodesic motion if and only if its covariant derivative $\nabla V$, as an endomorphism of $V^\perp$, has no real eigenvalues.
In particular, the set of unit vector fields defining geodesic motions is open in the $C^1$ topology. It follows that non-stationary geodesic motions do exist. 
\end{theorem}

\subsection{Hyperbolic manifolds}

\begin{theorem}\label{hypGeodMOTIONS}
A closed hyperbolic 3-manifold has no geodesic motion.
\end{theorem}

The proof (in Section \ref{hyperbolic}) consists in showing that a geodesic motion is stationary and then applying the non-existence result for geodesic foliations in this case \cite{ZeghibGeoFoliaofHypMan}. This is, however, not true for the hyperbolic space, which admits both stationary and non-stationary geodesic motions (see \ref{flexMotionOnHyp}). This phenomenon is true in fact for any Lie group with a left invariant metric (see section \ref{leftinvmet}).

\subsection{Codimension-one lightlike foliations}

In the constant curvature $1+2$-spaces, one has an interesting correspondence between one-dimensional null foliations and codimension-one geodesic lightlike foliations. 
\begin{theorem}\label{dimandcodim1_intro}
Let $\mathcal V$ be a smooth one-dimensional lightlike geodesic foliation of a Lorentzian $3$-manifold $(M,g)$. Then the orthogonal distribution $\mathcal V^\perp$ is integrable. If, in addition, $\mathcal V$ is complete and $(M,g)$ has constant sectional curvature, then the leaves of $\mathcal V^\perp$ are totally geodesic lightlike surfaces.
\end{theorem}

Therefore, the study of one-dimensional lightlike geodesic foliations of constant-curvature $(1+2)$-spacetimes is equivalent to the study of their foliations by complete geodesic lightlike surfaces. This observation confirms the result of Theorem \ref{thm_euclidean_motion_intro}, since, in general, foliations of $\R^n$ by hyperplanes are rigid, i.e., all leaves are parallel. In the positive-curvature case i.e. the de-Sitter space $\dS^{1,n-1}\subset \R^{1,n}$, any two lightlike hyperplanes must intersect in $\dS^{1,n-1}$, and thus no such foliation exist.

\subsection{Anti-de Sitter case}

The anti-de Sitter space $\AdS^{1,n}$ is the $-1$-sphere of the signature $(2,n)$ form on $\R^{2+n}$, its quotient under $\pm \Id$ defines the projective anti-de Sitter space $\P\AdS^{1,n}$. we have the following characterization:
\begin{theorem}\label{AdScodim1}
       Codimension-one lightlike geodesic foliations of $\widetilde{\AdS^{1,n}}$ are in one-to-one correspondence with inextendible causal curves in the universal Einstein universe $(\S^{n-1}\times\R, g_{o}-dt^2)$. 
\end{theorem}

We can further drop the $C^1$ regularity assumption in Theorem \ref{dimandcodim1_intro} if the geodesic lightlike foliation is global. More precisely,

\begin{theorem}
    Let $\mathcal{L}$ be a $C^0$ null foliation of $\widetilde{\AdS^{1,2}}$ (i.e. a $C^0$-foliation by lightlike geodesics). Then the orthogonal distribution $\mathcal{L}^\perp$ is tangent to a lightlike geodesic codimension-one foliation.
\end{theorem}

\begin{corollary}
    One-dimensional null geodesic foliations of $\AdS^{1,2}$ are characterized by causal loops of degree $1$ in $(\S^1 \times \S^1, dx^2-dy^2)$ i.e. graphs of distance-decreasing maps from $\S^1$ to itself.
\end{corollary}

The anti-de Sitter space $\AdS^{1,2}$ and even its compact quotients (\cite{KulkarniRaymond, salein2000varietes}) have very flexible timelike geodesic foliations (see \cite{DancigerGueritaudKassel, sagman} for sufficient existence conditions). In fact, this case exhibits a behavior similar to the one for foliations of $\S^3$ by great circles in \cite{gluck1983great}. Drawing this parallel, we can obtain, as already noticed in \cite{DancigerGueritaudKassel}, but with a different approach similar to the one in \cite{gluck1983great}, the following 
\begin{theorem}[c.f. \cite{DancigerGueritaudKassel, sagman}]

   Given a geodesic foliation $\F$ of $\widetilde{\AdS^{1,2}}$, Let $U$ be a connected component of the region where the leaves of $\F$ are timelike. If $U$ is non-empty, then $\F_{\mid U}$ is determined by the graph of a distance-decreasing or distance-increasing map
$f:\H^2\longrightarrow \H^2$.
Moreover, if $f$ is a contraction (i.e. $\operatorname{Lip}(f)\leq k<1$) then $\F$ is a global foliation i.e. $U=\widetilde{\AdS^{1,2}}$.
\end{theorem}

It is important to note that not all distance-decreasing maps give rise to global timelike geodesic foliations. There are examples of mixed foliations of $\AdS^{1,2}$ with spacelike, lightlike and timelike leaves. There are even mixed foliations with only lightlike and timelike leaves. Interestingly, for even-dimensional anti-de Sitter, there is no geodesic timelike foliations, even locally \cite{Tholozan_timelike_foliation}, this comes from a general fact about local foliations of even-dimensional spheres by great circles (see \cite{Ferusfoliationsofspheres}) and can be generalized for any $\AdS^{p,q}$, with $p+q$ even. Spacelike foliations of $\AdS^{1,2}$ are flexible too, as any foliation of the hyperbolic plane can be lifted to a spacelike foliation of $\AdS^{1,2}$. The situation however, becomes very rigid for compact $\AdS^{1,2}$-manifolds (see \cite{geodfoliationsofcompact}).

\subsection{AI declaration}
No generative AI tool has been used in any stage of this work. All the mathematical arguments have been developed without the use of AI, and the manuscript has been
entirely human-written.

\section{Further comments and results}
\subsection{Lie groups with a left invariant metric}
A left invariant metric on a Lie group is a Riemannian metric for which any left translation is an isometry. Such metrics are defined entirely by their restriction on the Lie algebra. This homogeneity gives a particular set of geodesic motions 

\begin{proposition}
    Let $(G,q)$ be a Lie group endowed with a left-invariant Riemannian metric $q$ and let $V$ be a unit left-invariant vector field. Then $V$ defines a geodesic motion. Moreover, there always exists a stationary geodesic motion.
\end{proposition}

If $\Gamma$ a lattice in $G$, $V\in \mathcal{X}(G)$ a left invariant vector field. By left invariance, both the metric and $V$ descend to the quotient $\Gamma \backslash G$ and hence $V$ defines a geodesic motion. 
In particular, for $G=\PSL(2,\R)$, $\Gamma$ cocompact, and $V$ the left invariant vector field defined by $\begin{pmatrix}
0 & 1\\
0 & 0
\end{pmatrix}$ at $\Id$. The flow of $V$ is called the horocyclic flow on $\Gamma \backslash \PSL(2,\R) \cong T^1\Sigma$, where $\Sigma= \Gamma \backslash \H^2$. $V$ is known (see \cite{sullivan1978foliation} or \cite[pp.~22--24]{Rechtman2009}) to be non-geodesible i.e. there is no Riemannian metric on  $T^1\Sigma$ for which the trajectories of $V$ are geodesic (even unparametrized). Notice however that by the proposition above $V$ defines a geodesic motion for any left invariant metric on $\PSL(2,\R)$.

\subsection{Laminations of hyperbolic spaces}
Let $q$ be the usual Lorentz form of $\R^{1,n}$. The $-1$-level i.e. $\{q(u)= -1\}$, is the hyperbolic space $\H^n$ and the $+1$-level is the de Sitter space $\dS^n$. A timelike hyperplane $H$ in $\R^{1,n}$ defines a geodesic hypersurface of $\H^n$. It orthogonal
cuts $\dS^n$ in two points, exactly one if one instead considers the projective version $\P\dS^n$. This is causal but non chronologically orientable. A codimension 1 geodesic lamination $\F$ of $\H^n$ is a family of non-intersecting geodesic hypersurfaces. It turns out that its corresponding subset in $\P\dS^n$ is a totally ordered subset, and conversely, any such totally ordered subset gives rise to a geodesic lamination.   It is however, more complicated to interpret the lamination property by means of the ``linear'' $\dS^n$. 
Finally, Theorem \ref{AdScodim1} extends to laminations of the universal anti-de Sitter space.  They correspond to totally ordered subsets of the universal Einstein universe.

\subsection{Codimension-one lightlike foliations of compact $\AdS$-manifolds}
\begin{figure}[h]
\includegraphics[scale=0.45]{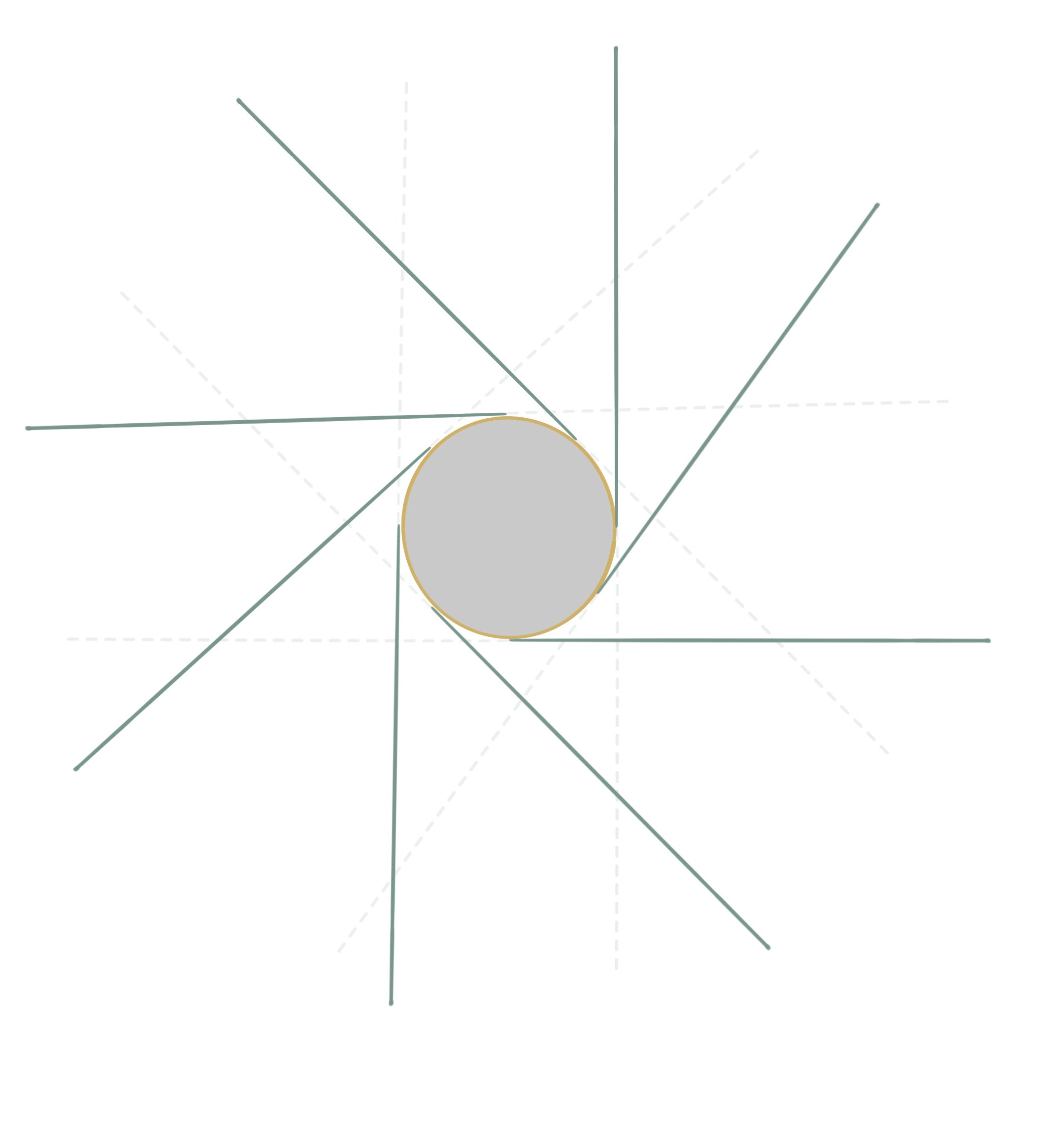} \centering
\caption{} \label{Foliation_of_Cylinder}
\end{figure}

An erroneous argument in the literature (\cite{zeghibAdS, saleinthese}) asserts that lightlike totally geodesic foliations of $\AdS^{1,n}$ and its compact quotients are rigid. The argument runs as follows. A lightlike  totally geodesic hypersurface is identified with the intersection $H_u=\AdS^{1,n}\cap u^\perp$, where $u$ is a nonzero null vector in $\mathbb{R}^{2,n}$. If $H_u$ and $H_v$ are disjoint, then
$ u^\perp\cap v^\perp
 =(\mathbb{R}u\oplus\mathbb{R}v)^\perp$
 does not meet $\AdS^{1,n}$ and therefore contains no timelike vectors. Consequently, the plane $\mathbb{R}u\oplus\mathbb{R}v$ must be totally isotropic. This leads to the claimed rigidity: the leaves are obtained by taking orthogonal complements of null lines in a fixed totally isotropic two-plane. \newline
The flaw is that $H_u$ has two antipodal connected components. A leaf is a single connected component, and two leaves may be disjoint even though the corresponding full intersections $H_u$ and $H_v$ meet. Thus, the argument establishes rigidity in the projective anti-de Sitter space, but not in $\mathrm{AdS}$ itself. \newline
A simple analogy illustrates the distinction (see Figure \ref{Foliation_of_Cylinder}). Let $C\subset\mathbb{R}^2$ be the unit circle, bounding the closed unit disk $\overline{B}$. The tangent lines to  $C$ do not foliate $\mathbb{R}^2\setminus\overline{B}$. However, after choosing an orientation of $C$, the open tangent rays pointing in the positive direction do form a foliation of this region. In section \ref{nullAdSsection}, we prove flexibility of codimension-one geodesic lightlike foliations of $\AdS^{1,n}$ and its covers. For compact quotients, we still have rigidity, the argument will be given in \cite{geodfoliationsofcompact}.

\subsection{Further remarks on null foliations in 4-dimensional Minkowski space}
Null geodesic foliations are called geodesic congruences in physical literature. They
play an important role in the Kerr-Schild construction. On a domain in four-dimensional Minkowski space,
consider a metric of the form
$
g=\eta+H\,\omega\otimes\omega,
$, 
where $\eta$ is the Minkowski metric, $H$ is a smooth function,
and $\omega=\eta(V,\cdot)$ is the metric dual of a nowhere-vanishing
null vector field $V$. We focus on the case where $V$ generates
a null geodesic foliation. Its integral curves are then also null
geodesics, up to reparameterization, for $g$. Once the normalization
of $V$ is fixed, the metric is determined by the foliation and
the function $H$. \newline
 Both the Kerr metric and its Schwarzschild limit admit such
a representation \cite{NewmanAdamo2014}. Their associated null
congruences are shear-free but expanding. Shear-freeness means
that the flow of $V$ preserves the conformal class of the metric
induced on the screen bundle $V^\perp/\langle V\rangle$.
The Schwarzschild congruence is twist-free, whereas the rotating
Kerr congruence has nonzero twist at generic points. Wherever
the twist is nonzero, the distribution $V^\perp$ fails the
Frobenius integrability condition. This distinguishes the
four-dimensional situation from dimension $1+2$, where null
geodesic congruences necessarily have vanishing twist. \newline
These particular congruences do not extend to smooth null
foliations of all of Minkowski space: the Schwarzschild
congruence has a caustic at the spatial origin, while the Kerr
congruence has a ring caustic. \newline
Shear-free null geodesic congruences also provide a fundamental
link with twistor theory. The Kerr theorem gives a local
description of real-analytic shear-free null geodesic congruences
in Minkowski space in terms of holomorphic surfaces in projective
twistor space, subject to appropriate regularity conditions
\cite{AdamoNewmanKozameh2012}. The term ``twistor,'' coined by
Penrose, itself reflects the twisting geometry of the shear-free
null geodesic congruences known as Robinson congruences
\cite{Penrose2016}.

\section{Geodesic motions in the Euclidean plane}

Let $U\subset \R^2$ be an open subset and $V:U\to \S^1$ be a unitary vector field on $U$. Define the map $F_t: p\mapsto p+ tV(p)$ for all $t\in \R$. We have 
\begin{theorem}\label{thm_euclidean_motion}

If $F_t$ is injective for all $t$, then the integral curves of $V^\perp$ are straight segments. In particular, when $U=\mathbb{R}^2$, the non-collision property i.e. $F_t$ is a geodesic motion, is equivalent to $V$ being constant.  

\end{theorem}
\begin{proof}[Proof of Theorem \ref{thm_euclidean_motion}]

For any $t\in\R$, the map $F_t$ is given by $$F_t: p\in U \mapsto p + t V(p) \in \R^2.$$ The non-collision condition implies that the map $F_t$ is a local $C^1$-diffeomorphism for any $t\in \R$. Hence for any $X=(x,y)\in\R^2$, the differential 
\[
D_XF_t=\Id+tD_X V
\]
is invertible for any $t\in\R$. We deduce that all the eigenvalues of $D_X V$ are equal to zero, hence $D_X V$ is nilpotent. 

After setting $V(X)=:(\cos(\theta(X)),\sin(\theta(X))$, we write 
\[
D_X V=\begin{pmatrix}
	 -\frac{\partial \theta }{\partial x}(X) \sin(\theta(X)) & 	 -\frac{\partial \theta }{\partial y}(X) \sin(\theta(X))  \\[4pt]
	 \frac{\partial \theta }{\partial x}(X) \cos(\theta(X)) & \frac{\partial \theta }{\partial y}(X) \cos(\theta(X)) 
	\end{pmatrix}.
\]
Observe that $\operatorname{Im}(D_X V)\subseteq \R\Bigl(-\sin(\theta(X)),\cos(\theta(X))\Bigr)=\R V(X)^\perp$.

Since $D_X V$ is nilpotent, we get $\operatorname{Im}(D_X V)\subset\ker(D_X V)$, that is, $V(X)^\perp\in \ker(D_X V)$. Hence, the vector field $V$ is constant on the integral curves of $V^\perp$. This implies that $V^\perp$ is constant along its own trajectories. Therefore, the integral curve of $V^\perp$ are straight segments.

Finally, when $U=\R^2$, the integral curves of $V^\perp$ define a foliation of $\R^2$ by affine lines. This could happen only if $V^\perp$ is constant. We conclude that the only collision-free $C^1$ vector field on $\R^2$ is the constant one.

\end{proof}

From Theorem \ref{thm_euclidean_motion}, one deduces that any global $C^1$ lightlike geodesic foliation of the 3-dimensional Minkowski space is given by parallel lines. Note that this is particular to the low dimensional case, as one can give many examples of non-constant lightlike foliations of $\R^{1,3}$. 
\begin{example}
    Consider the foliation of $\R^{1,3}$ by $\{x\}\times\R^{1,2}$ copies and fix on each copy a lightlike foliation given by a choice of a direction $V(x)\in\S^1$. Hence any smooth map $V:\R\rightarrow \S^1$ gives rise to a smooth null foliation of $\R^{1,3}$.
\end{example}

The aforementioned rigidity result for global null foliations of $\R^{1,2}$ is for $C^1$ lightlike foliations. The same conclusion holds in the $C^0$ case.

\begin{theorem}\label{C0rigidmink}
 Any lightlike geodesic $C^0$-foliation of $\R^{1,2}$ is constant.
\end{theorem}
\begin{proof}
   Let $\F$ be a $C^0$ lightlike foliation of $\R^{1,2}$ and let $P\subseteq \R^{1,2}$ be a spacelike plane. Each leaf of $\F$ intersects $P$ exactly once. We denote by $V:P\rightarrow \S^1$ the direction function with respect to the foliation $\F$ that associates to each point $p\in P$, the direction of the orthogonal projection on $P$ of the lightlike geodesic of $\F$ that passes through $p$. \par 
   Suppose by contradiction that $V$is not constant. Since $V$ is continuous, then there exists a level set of $V$that separates $P$, that is, there exists $s\in \S^1$ such that $P\setminus V^{-1}(s)$ is disconnected. There exists then a closed connected component $C\subseteq V^{-1}(s)$ that separates $P$. 

   \begin{lemma}
       The closed subset $C$ is convex.
   \end{lemma}
   \begin{proof}
      Fix $a,b\in C$ and $p\in[a,b]$. Denote by $l_p$ the lightlike geodesic through $p$ and by $\mathcal{F}_C$ the collection of parallel lightlike lines through $C\subset P$. Let $d_p$ denote the projection of $l_p$ on the plane $P$ with respect to the direction $V(C)\in \S^1$. So, the line $l_p$ intersects $\mathcal{F}_C$ if and only if $d_p$ intersects $C$. We deduce that $d_p$ is a single point, i.e. $l_p$ is parallel to the lines through $\F_C$. Indeed, if $d_p$ is an affine line then $d_p$ either intersects $C$, which is impossible since $l_p$ does not intersect $\F_C$, or $d_p$ separates $C$ which is also impossible sine $C$ is assumed to be connected. Thus, $d_p$ is a point.  It follows that $V(p)=V(C)$. In other words, we have $[a,b]\subset V^{-1}(s)$. Since $C$ is a connected component of $V^{-1}(s)$, we deduce that $[a,b]\subset C$. Hence $C$ is convex.

   \end{proof}
   \begin{lemma}
       Any closed connected convex subset $C$ that separates $\R^2$ is a band.
   \end{lemma}
   \begin{proof}
       Let us first show that for any point $p\in C$, there exists a line that goes through $p$ and is contained in $C$. Consider the following map 
       \begin{eqnarray*}
            \operatorname{depth}:\S^1&\rightarrow&[0,+\infty] \\
            \eta \, \, \,  &\mapsto& d(\eta):=\sup\{t\in[0,+\infty[\, |\, [p,p+t\eta]\subseteq C\}
       \end{eqnarray*}
       As $C$ is convex, we get that for any pair of disjoint points $\eta_1,\eta_2\in \operatorname{depth}^{-1}(+\infty)$, that are not antipodal, the geodesic arc $[\eta_1,\eta_2]$ is also contained in $\operatorname{depth}^{-1}(+\infty)$.

       On the other hand, as $C$ separates the plane, the subset $\operatorname{depth}^{-1}(+\infty)$ separates $\S^1$. This implies that $\operatorname{depth}^{-1}(+\infty)$ consists of two antipodals points $\{-\eta,\eta\}$. We conclude then that the line $(p+t\eta)_{t\in\R}$ lies in $C$. Hence, we conclude that $C$ is foliated by complete parallel lines which finishes the proof of the lemma.
       \end{proof}
       
We now have all the ingredients to prove Theorem \ref{C0rigidmink}. Let $L\subset \partial C$ be one of the two boundary lines of the band $C$. The leaves of $\F$ along $L$ define a lightlike affine plane $H$. On the other hand, any lightlike affine line in $\R^{1,2}$ either intersect $H$ or is parallel to the lightlike lines foliating $H$. This shows that all the leaves of $\F$ must be parallel to each other. That is, $V$ is constant. Contradiction!

\end{proof}
\section{Geodesic motions on spheres}
Let $\S^n$ be a round sphere. When $n$ is even, the hairy ball Theorem, implies that any vector field on $\S^n$ is singular, in particular, no motion can be defined. The situation is much more flexible in the odd dimensional case, as odd dimensional spheres admit many geodesic foliations (\cite{gluck1983great, CahnGluckNuchi2018}) and hence many stationary geodesic motions. In fact we have more
\begin{theorem}
    A unit vector field $V$ on $\S^n$ defines a geodesic motion if and only if its covariant derivative $\nabla V$, as an endomorphism of $V^\perp$, has no real eigenvalues.
In particular, the set of unit vector fields defining geodesic motions is open in the $C^1$ topology. It follows that a generic geodesic motion on the sphere is non-stationary. 
\end{theorem}
\begin{proof}
Let $V$ be a unit vector field on $\S^n$ ($n$ is necessary odd), and define the motion $F_t(p)=\exp_p(tV(p))$. Consider $u\in T_p\S^n$ non zero vector and a curve $c:]-\epsilon,\epsilon[$ tangent to it at $0$, we write $u=u^\perp + \lambda V(p)$, where $u\in V(p)^\perp$. The map $(t,s) \mapsto F_t(c(s)$ is a geodesic variation by definition, and the vector field $J(t)=\left.\frac{\partial }{\partial s}F_t(c(s) \right|_{s=0}$ is a Jacobi field along the geodesic $(F_t(p))_{t\in\R}$ and since $\S^n$ has constant curvature $1$, the Jacobi equation gives that 
$$J(t)=P_t( \lambda V(p)+ \cos(t)u^\perp+\sin(t) A(u)) $$
where $P_t$ is the parallel transport along $(F_t(p))_{t\in\R}$, and $A: w\in T_p\S^n \mapsto \nabla_w V(p)$ is the covariant derivative with respect to the Levi-Civita connection. Notice that since $V$ has norm $1$, the image of $A$ lies in $V(p)^\perp$.

Suppose that $V$ defines a geodesic motion i.e. $F_t$ is a diffeomorphism for all $t$, then, for any $u\in V(p)^\perp$ i.e. $\lambda=$, we have $J(t)\neq 0$ for any $t$ since $J(t)=d_pF_t(u)$, and this implies that $A: V(p)^\perp \rightarrow V(p)^\perp$ has no real eigenvalues. Indeed, if $J(t_0)=0$ for some $t_0$, then $\cos(t_0)u=-\sin(t_0)A(u)$, moreover, since $u\neq 0$, we get that $\sin(t_0)\neq 0$ and $-\frac{\cos(t_0)}{\sin(t_0}$ is an eigenvalue of $A$. Reciprocally, suppose the restriction of $A$ on $V(p)^\perp$ has no real eigenvalue, and consider $u\in T_p\S^n$, such that $J(t_0)=0$ for some $t_0$. Since $A(u)\perp V(p)$, we get 
$$P_{t_0}( \cos(t_0)u^\perp+\sin(t_0) A(u))=0 \;\; \text{and  } P_{t_0}(\lambda V(p))=0$$ 
which implies that $u\in V(p)^\perp$ and that it is an eigenvector of $A$ for the eigenvalue $-\frac{\cos(t_0)}{\sin(t_0}$, hence $u=0$, and $d_pF_t$ is invertible. We conclude that $F_t$ is a local diffeomorphism of $\S^n$, and since $\S^n$ is compact and simply connected, we get that $F_t$ is a diffeomorphism for all $t\in \R$.

Asking $A: w\in V^\perp \mapsto \nabla_w V^\perp$ to not have real eigenvalues is an open condition for the $C^1$ topology on $\mathcal{X}(\S^n)$, this comes from the fact that at each point $A$ depends only on the value of $V$ at its derivative at that point (this is a general property of the Levi-Civita connection, see \cite{OneilBook}). At the same time, a geodesic motion $V$ is stationary if $V$ is geodesic i.e. $\nabla_V V=0$, and this is a closed condition. Therefore, we deduce the existence of non-stationary geodesic motions on odd-dimensional spheres.

\end{proof}

\section{Geodesic motions on hyperbolic manifolds}\label{hyperbolic}

  Let $(M,\ip{}{})$ be a closed Riemannian three-manifold of constant sectional curvature  $-1$. 

\begin{theorem} There is no geodesic motion on $M$.

\end{theorem}

Assume by contradiction that there exists a collision-free unit vector field $V$ on $M$. We note the covariant derivative endomorphism  \(
 A : V^\perp\longrightarrow V^\perp,  A(Z)=\nabla_ZV
\) (it does preserve \( V^\perp\), since \(V\) is unit).

We first show that such a geodesic motion must be stationary, that is, $V$ is a geodesic flow (its orbits are geodesics). Then one uses \cite{ZeghibGeoFoliaofHypMan} to deduce that this is not possible. 
 
 \subsection{A volume growth argument}  We will have occasion to come back  in details 
 to this   following argument in \cite{geodfoliationsofcompact}, when we consider spacelike one-dimensional geodesic foliations of anti-de Sitter 3-manifolds. 
Consider the geodesic motion \( p \to F_t(p) = \exp_p t V(p) \).  Its derivative in the direction \( u \in V_p^\perp\) equals the $t$-evaluation of the Jacobi field $J$
along the geodesic $\gamma$ determined by \(V\) such that \(J(0) = u, \dot{J}(0)= A(u) \). Since $M$ has curvature $-1$, we get $$J(t) = \cosh t P_t(J(0)) + \sinh t P_t(\dot{J}(0))= \cosh t P_t(u) + \sinh t P_t(A(u)) $$ where $P_t$ is the parallel transport along the geodesic $F_t(p)$. As for the direction \(V_p\) itself, it gives the Jacobi field with \(J(0) = V_p\) and \( \dot{J}(0) = a = \nabla_VV\), which is orthogonal to \(V\); hence \( \ip{J(t)}{\dot \gamma(t)} \) is affine with initial values \(1, 0\), that is \( \ip{J(t)}{\dot \gamma(t)} = 1\). It follows, in a parallel frame adapted to \( \dot \gamma\), that the Jacobian of $F_t$ at $p$ is given by 
the determinant of  \( \cosh t \Id + \sinh t A_p \) acting on \( V_p^\perp\).

\begin{proposition}\label{Volume_preservation} Up to replacing $M$ by a finite cover, 
there are \(X,Y\)  
unit vector fields satisfying
\begin{equation}\label{eq:hyp}
 \ip{X}{V}=\ip{Y}{V}=0,\qquad
 \nabla_XV=X,\qquad \nabla_YV=-Y.
\end{equation}

\end{proposition}

 \begin{proof} Let \( \tau(p) = \tr A_p\) and \( \delta(p) = \det A_p\). The Jacobian of $F_t$ at $p$ is then 
 \[
 \det( \cosh t \Id + \sinh t A_p) = \cosh^2 t + \tau(p) \cosh t \sinh t + \delta (p) \sinh^2 t = \cosh^2 t \; \varphi_p (\tanh t),
 \]
 where \( \varphi_p(u) = \delta(p) u^2 + \tau(p) u + 1\). Being the Jacobian of a diffeomorphism, and equal to $1$ for $t = 0$, it is positive for any $t$. Therefore \(\varphi_p > 0\) on \( (-1,1)\), and so, by continuity, 
 \[
 \varphi_p(1) = 1 + \tau(p) + \delta(p) \geq 0, \qquad \varphi_p(-1) = 1 - \tau(p) + \delta(p) \geq 0.
 \]
 On the other hand, \(F_t\) being a global diffeomorphism, and $M$ having a finite volume, 
 \[
 \int_M \Jac F_t \; dv =  \mathsf{vol}(M) < \infty, \quad \text{for any } t \in \R.
 \]
 Now, \( \Jac F_t \sim \frac{1}{4} e^{2t} (1+ \tau + \delta) \) when \( t \to + \infty\). Would \( 1+ \tau + \delta\) be positive on a set of positive measure, Fatou's lemma would give \( \int_M \Jac F_t \; dv \to \infty\). Hence \( 1+ \tau + \delta = 0 \) almost everywhere, and everywhere by continuity. Letting \( t \to - \infty\) yields in the same way \( 1 - \tau + \delta = 0\).
 
 \; Therefore \( \tau = 0\) and \( \delta = -1\): the characteristic polynomial of $A_p$ is \( \lambda^2 -1\), and so $A_p$ has exactly $+1$ and $-1$ as eigenvalues. (In particular, the case where $A_p$ has a non-real eigenvalue does not occur, and \( \Jac F_t = 1\): the $F_t$ are volume preserving.) The eigenvalues being distinct and constant, $A$ is diagonalizable and its two eigenline fields are smooth; up to replacing $M$ by a finite cover, they are generated by unit vector fields $X$ and $Y$ satisfying \eqref{eq:hyp}.
 
 \end{proof}

\bigskip

Set the vector field $a = \nabla_V V$. We want to prove $a = 0$. Recall that the Riemann curvature tensor is defined in general by:
\[
 R(U,W)Z=\nabla_U\nabla_WZ-\nabla_W\nabla_UZ-\nabla_{[U,W]}Z.
\]
In the constant sectional curvature \(-1\), it has the special form
\[
 R(U,W)Z=-\bigl(\ip{W}{Z}U-\ip{U}{Z}W\bigr).
\]

\subsubsection{A first preparing case}  
Here, we want to prove   \(a = \nabla_VV=0\), assuming  $\ip{X}{Y} = 0$

Write \(a =\nabla_VV=kX+\ell Y\).
By orthonormality, there are functions \(\alpha,\beta\) such that
\[
 \nabla_XY=\alpha X,\qquad \nabla_YX=\beta Y.
\]
The curvature identity \(R(X,Y)V=0\) gives
\[
 0=-\alpha X-\beta Y-\nabla_{\alpha X-\beta Y}V
   =-2\alpha X-2\beta Y.
\]
Hence
\[
 \nabla_XY=\nabla_YX=0,\qquad
 \nabla_XX=-V,\qquad \nabla_YY=V.
\]

Write \(\nabla_VX=qY-kV\); then
\(\nabla_VY=-qX-\ell V\), and
\[
 [X,V]=X-qY+kV,\qquad [Y,V]=-Y+qX+\ell V.
\]
Taking respectively the \(X\)-component of \(R(X,V)V=-X\)
and the \(Y\)-component of \(R(Y,V)V=-Y\), we obtain
\begin{align*}
 -1&=\ip{R(X,V)V}{X}=X(k)-1-k^2,\\
 -1&=\ip{R(Y,V)V}{Y}=Y(\ell)-1-\ell^2.
\end{align*}
Thus
\[
  X(k)=k^2,\qquad Y(\ell)=\ell^2
\]

 By compactness, \(k^2\) attains its maximum at some point \(p\).
At that point,
\[
 0=X(k^2)(p)=2k(p)X(k)(p)=2k(p)^3.
\]
Therefore the maximum of \(k^2\) is zero, so \(k\equiv0\).
The same argument applied to \(\ell^2\), using \(Y(\ell)=\ell^2\),
gives \(\ell\equiv0\). Consequently \(\nabla_VV=0\).

\subsubsection{The general case} Now, we no longer assume $X$ is orthogonal to $Y$.   We replace 
$X$ by the new vector field $N$, such that 
\begin{equation}\label{eq:frame}
 c=\ip{X}{Y},\qquad s=\sqrt{1-c^2}>0,\qquad
 N=\frac{X-cY}{s},\qquad h=\frac{2c}{s}.
\end{equation}
(The inequality \(s > 0\) holds since $X$ and $Y$, being eigenvectors of $A$ for distinct eigenvalues, are independent.) Then \((Y,N,V)\) is an orthonormal frame and \(X=cY+sN\). As above, 
write
$
 a=\ell Y+kN
$. 
(So the function \(k=\ip{a}{N}\) is the component of the acceleration
perpendicular to \(Y\)).

\medskip

Like the previous starting case, we have 
\begin{lemma}\label{lem:identity} We have: 
 ${X(k)=s\,k^2}$.
\end{lemma}

\begin{proof}

Consider the differential form defined by   \(\omega(Z)=\ip{\nabla_ZY}{N}\).
 
From definitions and  \(X=cY+sN\), we infer 
\[
 \nabla_YV=-Y,\qquad \nabla_NV=N+hY,\qquad
 \nabla_VV=\ell Y+kN.
\]
It follows that
\begin{align*}
 \nabla_ZY&=\omega(Z)N-\ip{Y}{\nabla_ZV}V,\\
 \nabla_ZN&=-\omega(Z)Y-\ip{N}{\nabla_ZV}V.
\end{align*}
In particular,
\begin{align*}
 \nabla_YY&=\omega(Y)N+V,
 &\nabla_YN&=-\omega(Y)Y,\\
 \nabla_NY&=\omega(N)N-hV,
 &\nabla_NN&=-\omega(N)Y-V,\\
 \nabla_VY&=\omega(V)N-\ell V,
 &\nabla_VN&=-\omega(V)Y-kV.
\end{align*}

We now use the three curvature identities
\[
 R(Y,N)V=0,\qquad R(Y,V)V=-Y,\qquad R(N,V)V=-N.
\]
Taking their \(N\)-components, calculated from the previous connection formulas, 
\begin{align*}
 \ip{R(Y,N)V}{N}
 &=h\omega(Y)+2\omega(N)-hk,\\
 \ip{R(Y,V)V}{N}
 &=Y(k)+\ell\omega(Y)+2\omega(V)-\ell k,\\
 \ip{R(N,V)V}{N}
 &=N(k)+\ell\omega(N)-h\omega(V)-1-k^2.
\end{align*}
Equating these expressions to \(0,0,-1\), respectively, yields
\begin{align}
 2\omega(N)&=h\bigl(k-\omega(Y)\bigr),\label{eq:omega}\\
 Y(k)&=\ell k-\ell\omega(Y)-2\omega(V),\label{eq:Yk}\\
 N(k)&=k^2-\ell\omega(N)+h\omega(V).\label{eq:Nk}
\end{align}
Since \(sh=2c\), equation \eqref{eq:omega} implies
\[
 \omega(X)=c\omega(Y)+s\omega(N)=ck.
\]
Combining \eqref{eq:Yk} and \eqref{eq:Nk}, we obtain
\begin{align*}
 X(k)
 &=cY(k)+sN(k)\\
 &=sk^2+\ell\bigl(ck-c\omega(Y)-s\omega(N)\bigr)
     +(sh-2c)\omega(V)\\
 &=sk^2.
\end{align*}
 
\end{proof}

Exactly as in the preparing case, $k$ vanishes, which means that $a = \nabla_ V V$ is collinear to $Y$. Now, $-V$ is a geodesic motion as well, since \(F_t^{-V} = F_{-t}\), with the same acceleration \( \nabla_{-V}(-V) = a\), and with endomorphism $-A$: the roles of $X$ and $Y$ are thus exchanged. Applying what precedes to $-V$, that is, replacing $N$ by \( N' = \frac{Y - cX}{s}\), yields that $a$ is collinear to $X$. Since $s >0$, $X$ and $Y$ are independent, and hence $a = 0$.
$\Box$

\begin{example}[Non-stationary motion on $\H^2$]\label{flexMotionOnHyp} Consider the upper half plane $\H^2$ with the hyperbolic metric. Then $\Aff^+(\R)\subset\PSL(2,\R)$ acts freely transitively on $\H^2$. Now fix $V_0=\partial x$ at the point $i\in\H^2$ and translate it using elements of $\Aff^+(\R)$. This gives a vector field $V\in\mathcal{X}(\H^2)$ with integral curves the horizontal lines in $\H^2$. We claim that $V$ defines a geodesic motion on $\H^2$. Indeed, consider the geodesic $\gamma$ tangent to $V_0$ at $i$. The motion $(\exp(tV(p)),t)=(g\gamma(t),t)$ where $g$ is the unique element of $\Aff^+(\R)$ sending $i$ to $p$. One can easily prove that this defines a lightlike foliation of $\H^2 \times \R$. For instance, if $(g\gamma(t),t)=(h\gamma(s),s)$, then $t=s$ and $g=h$ by freeness of the action. In particular, $V$ defines a non-stationary geodesic motion on $\H^2$.
\end{example}

 \section{Geodesic motions on Lie groups}\label{leftinvmet}
 
Example \ref{flexMotionOnHyp} can be generalized to the more general context of Lie groups. Recall that a left-invariant metric on a Lie group is a metric invariant under left translation. We have

\begin{proposition}
    Let $(G,q)$ be a Lie group endowed with a left-invariant Riemannian metric $q$ and let $V$ be a unit left-invariant vector field. Then $V$ defines a geodesic motion. Moreover, there always exists such a vector field generating a stationary geodesic motion.
\end{proposition}
\begin{proof}
   First, $(G,q)$ is homogeneous and, in particular, geodesically complete. Let $V$ be a unit left-invariant vector field. For $g\in G$, let $\gamma_g(t)$ be the geodesic starting at $g$ with speed $V(g)$. We have $\gamma_g(t)=g.\gamma_e(t)$ where $e$ denotes the identity element of $G$. It follows that $\gamma_g(t)=\gamma_{g^\prime}(t)$ if and only if $g=g^\prime$. In other words, $V$ is collision-free.

   Moreover, there always exists a non-zero (unit) vector $v\in \g$ such that the one-parameter subgroup $\exp(tv)$ is a geodesic for $q$ (see e.g. \cite[Prop. 1.]{Bromberg_2008}). Thus, the left-invariant vector field generated by $v$ has geodesic orbits, i.e., it generates a stationary geodesic motion.
\end{proof}

\begin{remark}
    Consider the universal cover of the Euclid group $\widetilde{\Isom^+}(\R^2)$. It can be identified with the semi-direct product $\R^2 \rtimes \R$, where $\R$ acts on $\R^2$ by rotation. The usual Euclidean metric on $\R^3$ is a left-invariant metric on  $\widetilde{\Isom^+}(\R^2)$. By the proposition above,  the following vector field \[
 V(x,y,z)=\frac1{\sqrt2}
 \left(\cos(2\pi z), \sin(2\pi z), 1\right).
\] 
is left invariant and hence defines a geodesic motion on $\R^3$. Moreover, this vector field descends to the quotient $\T^3=\R^3/\Z^3$
\end{remark}

\color{black}
\section{Null foliations in constant curvature $(1+2)$-spaces}
Let $(M,g)$ be a Lorentzian three-manifold and consider $\mathcal V$ a foliation of $M$
by lightlike geodesics generated by a vector field $V$ such that $\nabla_VV=0$ .

\begin{proposition}
    The orthogonal distribution $V^\perp$ is integrable.
\end{proposition}

\begin{proof}
Let $X$ be a non-singular spacelike vector field such that $V^\perp=\operatorname{span}\{V,X\}$. By Frobenuis' theorem and because $\dim(M)=3$, we obtain that $ V^\perp$ is integrable if and only if $[X,V] \in V^\perp$, i.e. $g([X,V],V)=0$.
We have \begin{align*}
    g([X,V],V) =& g(\nabla_X V-\nabla_V X,V)\\
    =& g(\nabla_X V,V) -g(\nabla_V X,V) \\
    =&\frac{1}{2} X.g(V,V) - V.g(X,V) +g(X,\nabla_V V) \\
    =& 0.
\end{align*}
Where we used only the facts that $g(V,V)=g(X,V)=0$ and $\nabla_V V=0$. 
\end{proof}

\begin{remark} There is in fact an equivalence between one-dimensional null foliations with integrable orthogonal distributions and codimension-one lightlike foliations. Examples of the latter foliations are given by the level sets of a solution of the eikonal equation (see for instance \cite{GourgoulhonJaramillo2006, Mounoud2004}).
\end{remark} 
\begin{theorem}\label{Dim1Codim1}
       Suppose that $(M,g)$ has constant sectional curvature and that $V$ is complete. Then the lightlike foliation $\mathcal V^\perp$ is totally geodesic.
\end{theorem}

 \begin{proof}
    Let $S$ be a leaf of $\mathcal V^\perp$. Since $V\in V^\perp$, $S$ is foliated by the complete lightlike geodesics tangent to $V$. Fix $p\in S$ and choose a smooth curve $c:(-\varepsilon,\varepsilon)\longrightarrow S$ such that
$$ c(0)=p, \qquad c'(s)\notin\mathbb{R}V(p), \;\; \forall s\in (-\varepsilon,\varepsilon). $$
Notice in particular that $c'(s) \perp V$. Let $\varphi_t$ be the flow of $V$ and define a parametrization $F(s,t)=\varphi_t(c(s)).$ For every fixed $s$, the curve
$ t\longmapsto F(s,t) $
is a complete parametrized lightlike geodesic. Consider the geodesic $\gamma(t)=F(0,t)$ and define the variational vector field
$$ J(t)=\left.\frac{\partial F}{\partial s}\right|_{s=0}. $$
Since $F$ is a geodesic variation of $\gamma$, we have that $J$ is a Jacobi field along $\gamma$ with $J(0)=c'(0)$. Moreover, $J(t)\in V(t)^\perp$ (since it is tangent to $S$) and $J(t)$ is never collinear with $V(t)$, since $F$ is a local parametrization of the leaf $S$.

Recall that $J$ satisfies the Jacobi equation:
$$ \nabla_V^2J+R(J,V)V=0.$$

Since $M$ has constant sectional curvature $\kappa$, we have
$$ R(A,B)C=\kappa\bigl( g(B,C) A - g( A,C) B \bigr).$$

Taking $A=J$, $B=C=V$, and using
$$ g(V,V)=0, \qquad g(J,V)=0,$$
we get $R(J,V)V=0$. and hence the Jacobi equation gives
$$\nabla_V^2J=0.$$

Thus $\nabla_VJ$ is parallel along $\gamma$. If $P_t$ denotes parallel transport along $\gamma$, then
$$J(t)=P_t\left(J(0)+t\nabla_VJ(0)\right).$$

Choose a parallel spacelike vector field $X$ along $\gamma$ such that $V^\perp=\operatorname{span}(V,X)$ and so that $X(0)=J(0)$. In the basis $\{V,X\}$ of $V^\perp$ we can write:
$$\nabla_VJ(0)=aV(0)+bX(0).$$
Since $V$ and $X$ are parallel along $\gamma$, the transverse
component of $J(t)$ is $(1+bt)X(t)$. If $b\neq0$, then for $t_0=-\frac1b\in\mathbb{R}$ we get $J(t_0)\in\mathbb{R}V(t_0)$. This contradicts the fact that the variation $F$ parametrizes the
surface $S$ and hence that $J(t)$ is everywhere transverse to $V(t)$. Therefore $b=0$, hence
$$\nabla_VJ(0)\in\mathbb{R}V(0).$$
Since $\nabla_VJ$ is parallel along $\gamma$ and $V$ is geodesic, it follows that
\begin{equation}\label{equa1}
 \nabla_VJ\in\mathbb{R}V  
\end{equation}
everywhere along $\gamma$.

Now we have, $[V,J]=0$ because $V$ and $J$ arise as the coordinate vector fields of the parametrization $F$. Since the Levi-Civita connection is torsion-free, we get $\nabla_JV=\nabla_VJ$. By equation \ref{equa1}, we get
$\nabla_JV\in\mathbb{R}V.$
Using this and the fact that $\nabla_VV=0$, we obtain
$$\nabla_YV\in\mathbb{R}V \qquad \text{for every }Y\in TS.$$

Finally, let $Y,Z\in TS$. Since $TS=V^\perp$, we have $g(Z,V)=0.$
Differentiating this identity in the direction $Y$ gives
$$0=Y(g(Z,V))=g(\nabla_YZ,V)+g(Z,\nabla_YV).$$
Since $\nabla_YV$ is proportional to $V$ and $Z\perp V$, we have $g(Z,\nabla_YV)=0$. Therefore $g(\nabla_YZ,V)=0.$
Hence $\nabla_YZ\in V^\perp=TS$, which implies that $S$ is totally geodesic.
\end{proof}

In the positive curvature case, the situation is rigid. Recall that the de Sitter space $\dS^{1,n}$ is defined as the unit ball (radius $1$) in Minkowski space $\R^{1,n}$. A complete totally geodesic lightlike hypersurface is the intersection of $\dS^{1,n}$ with a degenerate hyperplane of $\R^{1,n}$. Notice that any two such hyperplanes must intersect in a spacelike subspace of codimension $2$ and, consequently, two lightlike geodesic hypersurfaces in $\dS^{1,n} $ must intersect. Therefore, $\dS^{1,n} $ does not admit a codimension-one lightlike geodesic foliation. In particular, by the theorem above, $\dS^{1,2}$ admits no one-dimensional null geodesic foliation.

The anti-de Sitter case is rich and more involved. The rest of the paper is therefore devoted to the study of the $\AdS$ case.

\section{Lightlike geodesic foliations of anti-de Sitter spaces}\label{nullAdSsection}

Let $q$ be a quadratic form on $\R^{2+n}$ of signature $(2,n)$. The anti-de Sitter spacetime $\AdS^{1,n}$, is defined as the $-1$-level of $q$ endowed with the induced (Lorentz) metric. The isometry group of $\AdS^{1,n}$ is therefore the algebraic group $\O(2,n)$. 

The group $\O(2,n)$ contains $-\Id$ as a central element. If follows that $\PO(2,n)$ is the full isometry group of the quotient $\P\AdS^{1,n}=\{\pm\Id\}\backslash\AdS^{1,n}$. We call $\P\AdS^{1,n}$ the \textit{projective} anti-de Sitter spacetime.

\subsection{The conformal boundary of \texorpdfstring{$\AdS^{1,n}$}{ADSn}}

In $\R^{2,n}$, consider the light cone $C$ of the quadratic form $q$ i.e. 
$q^{-1}(0)\setminus \{0_{\R^{2,n}}\}$. The flat metric defined by $q$ restricts to a degenerate metric on $C$ of signature $(0,-,+,..,+)$. Moreover, at any point $p\in C$ the degenerate direction of $q|_{T_pC}$ is the kernel of $d_p\pi$, where $\pi : \R^{2,n}\backslash \{0_{\R^{2,n}}\}\rightarrow \S^{n+1}$ is the natural projection. This gives for instance, a conformal structure on $\pi(C)$ as tangent spaces on each fiber are identified via homotheties. Furthermore, $\pi(C)$ is the boundary of the open domain $\pi(\AdS^{1,n})$ which is identified with $\AdS^{1,n}$, since the restriction of the projection $\pi$ on $\AdS^{1,n}$ is a smooth embedding. The space $\pi(C)$ thus merits the name of conformal boundary of $\AdS^{1,n}$, noted $\widehat{\Ein}^n$. If we use the projectivization to $\R\mathbb{P}^{n+1}$, then the restriction on the projection to $\AdS^{1,n}$ will give an open domain in $\R\mathbb{P}^{n+1}$, that identifies naturally with $\P\AdS^{1,n}$. Moreover, the light cone $C$, projects into a boundary of this open domain (hence to a boundary of $\P\AdS^{1,n}$), denoted by $\Ein^n$ and called the Einstein universe. Notice that $\widehat{\Ein}^n$ is a double cover of $\Ein^n$. 

On $\R^{2,n}$, fix the coordinates $(x_1,x_2, y_1,..,y_n)$, where $C$ is defined by the equation $-x_1^2-x_2^2+y_1^2+..+y_n^2=0$. The space $\widehat{\Ein}^n$ is naturally identified with the intersection in $\R^{2,n}$ of $C$ with the sphere of equation $x_1^2+x_2^2+y_1^2+..+y_n^2=2$. Thus, $\widehat{\Ein}^n$ can is defined by the equations $y_1^2+..+y_n^2=1$ and $x_1^2+x_2^2=1$ and consequently, $\widehat{\Ein}^n$ is identified with $\S^1\times \S^{n-1}$. Moreover, the restriction of $q$ in $\widehat{\Ein}^n$ endows it with the lorentzian metric $-dt^2+g_{\S^{n-1}}$ in the product $\S^1\times \S^{n-1}$ and which defines the same conformal structure as before (see \cite{Frances2002}).

For $n\geq3$, the universal cover of $\Ein^n$ (and of $\widehat{\Ein}^n$) denoted $\widetilde{\Ein}^n$ is identified with $\R \times \S^{n-1}$ with the conformal class $[-dt^2+g_{\S^{n-1}}]$. Considering the deck transformation $\tau: (t,x)\in \widetilde{\Ein}^n \cong \R \times \S^{n-1} \mapsto (t+\pi, -x)$, we get $\Ein^n=\widetilde{\Ein}^n/<\tau>$ and $\widehat{\Ein}^n=\widetilde{\Ein}^n/<\delta>$, where $\delta=\tau^2$. In the $3$-dimensional case, the universal cover is identified with $\R^2$, however, it is more helpful to consider the $\Z$-covering $\R\times\S^1$ denoted for simplicity $\widetilde{\Ein}^2$ and we get as before $\Ein^2=\widetilde{\Ein}^2/<\tau>$ and $\widehat{\Ein}^2=\widetilde{\Ein}^2/<\delta>$. 

The group of isometries of $\AdS^{1,n}$ is given by $\O(2,n)$. The same group acts on the cone $C$ isometrically, and therefore it acts conformally on $\widehat{\Ein}^n$ and on $\Ein^n$. The action on $\Ein^n$ becomes faithful if we quotient by $\pm \Id$, and it can be proven that $\P\O(2,n)$ is the full conformal group of $\Ein^n$ (see \cite{Frances2002}).

\subsubsection{\textbf{The conformal embedding}}

One can see the Einstein universe as the conformal boundary of $\AdS^{1,n}$ via a conformal embedding 
\[
\AdS^{1,n}\hookrightarrow\widehat{\Ein}^{n+1}\simeq\mathbb \S^1\times \S^{n}.
\]
Indeed, recall that the anti-de Sitter $\AdS^{1,n}$ is the $-1$ sphere: $q^{-1}(-1)$ for the $(2,n)$ form $q$ on $\R^{2,n}$ , now multiply $\R^{2,n}$ by a spacelike $(\R, dt^2)$, the map $x\in\AdS^{1,n}\hookrightarrow (x,1)\in \R^{2,n}\times \R$ embeds $\AdS^{1,n}$ isometrically in the isotropic cone $C^{n+2}$ of $\R^{2,n}\times \R$, composing with the projection to $\widehat{\Ein}^{n+1}$ gives the desired conformal embedding. Hence, $\AdS^{1,n}$ is identified with the open domain $U_+=\{[x:1]\in \widehat{\Ein}^{n+1} \}$ in the space of rays $\S^{n+2}$, and its boundary is exactly the image of $\widehat{\Ein}^n$ identified with the subset $\{ [x:0] \in \widehat{\Ein}^{n+1} \; q(x)=0\}$. After projectivization, we get a conformal embedding $\P\AdS^{1,n}\hookrightarrow\Ein^{n+1}$ with boundary $\Ein^n$ and lifting everything to universal cover (or the $\Z$-cover in dimension $1+n=3$) we get a conformal embedding $\widetilde{\AdS^{1,n}}\hookrightarrow\widetilde{\Ein}^{n+1}$ (see \cite{Frances2002, Francesconfboundary, Barbotsurvey} for more details).

We can see this conformal embedding more concretely, consider the unit ball $\B^{n}=\{(x_1,..,x_{n}), \; ||x||^2:=\sum_{i=1}^{n-1} x_i^2 < 1 \} $ endowed with the hyperbolic metric $g_{hyp}$. we can verify (see \cite{SeppiBonsante}) that the map $$(x,t)\in \S^1 \times \B^n \mapsto \frac{1}{1-||x||^2}\left(2x, (1+||x||^2)\cos(t), (1+||x||^2) \sin(t)\right) \in \AdS^{1,n}$$
is a global isometry, where $\S^1 \times \B^n$ is endowed with the Lorentzian metric $g_{lor}:=g_{hyp}- \left(\frac{1+||x||^2 }{1-||x||^2}\right)^2 dt^2 $. Now, notice that the metric $\left(\frac{1-||x||^2 }{1+||x||^2}\right)^2g_{hyp}$ is exactly the round metric on $\B^{n}$, seen as a hemisphere of $\S^{n}$. Up to lifting everything to the universal cover, we get
\begin{proposition}
The anti-de Sitter space $\widetilde{\AdS^{1,n}}$ is conformally equivalent to the open domain $
U=\mathbb R\times H_n,$
where $H_n$ is a hemisphere of $\S^{n}$, i.e. $\partial H_n$ is a great hypersphere $\S^{n-1}$.

\end{proposition} 
\begin{remark}
    Notice that taking the restriction of the covering map $\widetilde{\rho}:\widetilde{\Ein}^{n+1} \rightarrow \widehat{\Ein}^{n+1} $ to the image of $\widetilde{\AdS^{1,n}}$ under the conformal embedding above gives the cyclic cover $\widetilde{\AdS^{1,n}}\rightarrow \AdS^{1,n}$.
\end{remark}
 
\color{black}
\subsubsection{End points of lightlike geodesics}

Let $\P\mathcal L$ denote the space of maximal lightlike geodesics of $\P\AdS^{1,n}$. Recall that $
\Ein^{n}=\mathbb P(q^{-1}(0)),$ the projectivization of the null cone. There is a natural projection
\[
\Omega:\P\mathcal L\longrightarrow \Ein^{n}
\]
that associates to each lightlike geodesic $\ell\in\mathcal L$ the isotropic line that directs it. More precisely, recall that any plane $P$ of signature $(0,-)$ contains a unique lightlike vectorial line $\ell_0$ and it intersects $\AdS^{1,n}$ along two affine lines $\ell_1$, $\ell_{-1}$ parallel to $\ell_0$. Moreover, both $\ell_1$ and $\ell_{-1}$ are lightlike geodesic, and any lightlike geodesic of $\AdS^{1,n}$ is obtained in this way. Quotienting by $\pm \Id$ identifies the two geodesics $\ell_1$ and $\ell_{-1}$ and defines a unique lightlike geodesic $\mathbb{P}(P)$ in $\P\AdS^{1,n}$. Notice that this map is equivariant under the conformal group $\P\O(2,n)$.
\begin{definition}For a lightlike geodesic $\ell\in\mathcal L$, the element $\Omega(\ell)\in\Ein^{n}$ is called the \textit{endpoint} or the \textit{limit point} of $\ell$.

\end{definition}

Fix a time orientation on $\P\AdS^{1,n}$, and define the space $\mathcal{L}$ of future oriented lightlike geodesics of $\AdS^{1,n}$. We can define two maps $\widehat{\Omega}_+$, $\widehat{\Omega}_-$  on $\mathcal{L}$, which for each geodesic $\ell$ directed by an isotropic vector $v$,  associate to it its future end point $\widehat{\Omega}_+(\ell)=[v]\in \widehat{\Ein}^n$ and its past endpoint $\widehat{\Omega}_-(\ell)=[-v]\in \widehat{\Ein}^n$. If $\pi: \AdS^{1,n}\rightarrow\P\AdS^{1,n}$ and $\rho : \widehat{\Ein}^n \rightarrow \Ein^n$ are the double coverings, then we have for any $\ell \in \mathcal{L}$, $\rho(\widehat{\Omega}_+(\ell))=\rho(\widehat{\Omega}_-(\ell))=\Omega(\pi(\ell))$. Similarly, we can define the future and past endpoint maps $\widetilde{\Omega}_\pm$ from the space $\widetilde{\mathcal{L}}$ of future oriented lightlike geodesic of $\widetilde{\AdS^{1,n}}$ to $\widetilde{\Ein}^n$.

The conformal identification of $\widetilde{\AdS^{1,n}}$ with $\R^1 \times H_n$ endowed with the metric $-dt^2+g_o$ gives a nice description of the endpoint maps $\widetilde{\Omega}_\pm$. A general fact in Lorentzian geometry (see \cite[Appendix~D, p.~446]{Wald1984}) states that unparametrized lightlike geodesics are still geodesics after a conformal change of the metric. Thus, a future directed lightlike geodesic $\gamma$ in  $\widetilde{\AdS^{1,n}}$ is given by $t\mapsto (t,c(t))$ for $a<t<a+\pi$ and where $c$ is a geodesic in $H_n$. We get $\widetilde{\Omega}_-(\gamma)=(a,c(a))$ and $\widetilde{\Omega}_+(\gamma)=(a+\pi, -c(a))$. In particular, we have $\tau\circ\widetilde{\Omega}_-=\widetilde{\Omega}_+$.

\begin{remark}
    For a point $\ell_0\in\Ein^{n}$, the fiber $
\Omega^{-1}(\ell_0)=F_0$ admits a natural Riemannian metric isometric to the hyperbolic space $\mathbb H^{n-1}$. Indeed, $\Omega^{-1}(\ell_0)$ is the collection of all affine lines parallel to $\ell_0$ and contained in $\P\AdS^{1,n}$. In other words, the elements of $\Omega^{-1}(\ell_0)$ are the leaves of the canonical foliation of the codimension-one lightlike degenerate geodesic hypersurface $
\mathbb{P}(\ell_0^{\perp}\cap\AdS^{1,n}).$ 
This hypersurface admits a transverse Riemannian structure isometric to $\mathbb H^{n-1}$.

We get that $\mathcal L$ is a fiber bundle over the Einstein universe $\Ein^{n}$ with fibers isometric to the hyperbolic space $\mathbb H^{n-1}$. The group $\mathrm{PO}(2,n)$ preserves this structure, i.e. it acts isometrically in the vertical direction and conformally on $\Ein^{n}$. Observe that for $n=2$, a fiber of $\mathcal L$ is isometric to the Euclidean affine real line $\mathbb R$.

\end{remark}

\begin{fact}Connected totally geodesic complete lightlike hypersurfaces of $\widetilde{\AdS^{1,n}}$ are of the form
\[
\widetilde{\Omega}_-^{-1}(p),\qquad p\in\widetilde{\Ein}^{n}\simeq\mathbb R\times \S^{n-1}.
\]
Thus codimension-one totally geodesic lightlike hypersurfaces of $\widetilde{\AdS^{1,n}}$ are in one-to-one correspondence with points in the Einstein universe $\widetilde{\Ein}^{n}$.

\end{fact} 

With the identification $\widetilde{\Ein}^n\cong\R\times \S^{n-1}$, for any point $p=:(t,x)\in \widetilde{\Ein}^n$ we have :  
$$\widetilde{\Omega}_-^{-1}(p) = \{(t+s,y)\in \R\times H_n, \;\text{with} \; d(x,y)=s \}   $$

Moreover, we can characterize intersections of these hypersurfaces:
\begin{lemma}\label{intersection}
    For $p:=(t,x)$ and $q:=(t',x')$ in $\widetilde{\Ein}^n$ with $t\leq t'$, we have $\widetilde{\Omega}_-^{-1}(p)\cap \widetilde{\Omega}_-^{-1}(q)$ is empty, if and only if  $\mid t-t'\mid >d(x,x')$ i.e. $p$ and $q$ are causally related.
\end{lemma}
\begin{proof}
    An intersection point of $\widetilde{\Omega}_-^{-1}(p)$ and $\widetilde{\Omega}_-^{-1}(q)$ is given by a point $y\in H_n$, that satisfies
$$t+d(x,y)=t'+ d(x',y)$$
The triangle inequality implies therefore that  $\mid t-t'\mid  \leq d(x,x') $, and if the latter inequality is true then we can find $y$ on the geodesic between $x$ and $x'$.
\end{proof}

\begin{figure}[h] 
\includegraphics[scale=0.5]{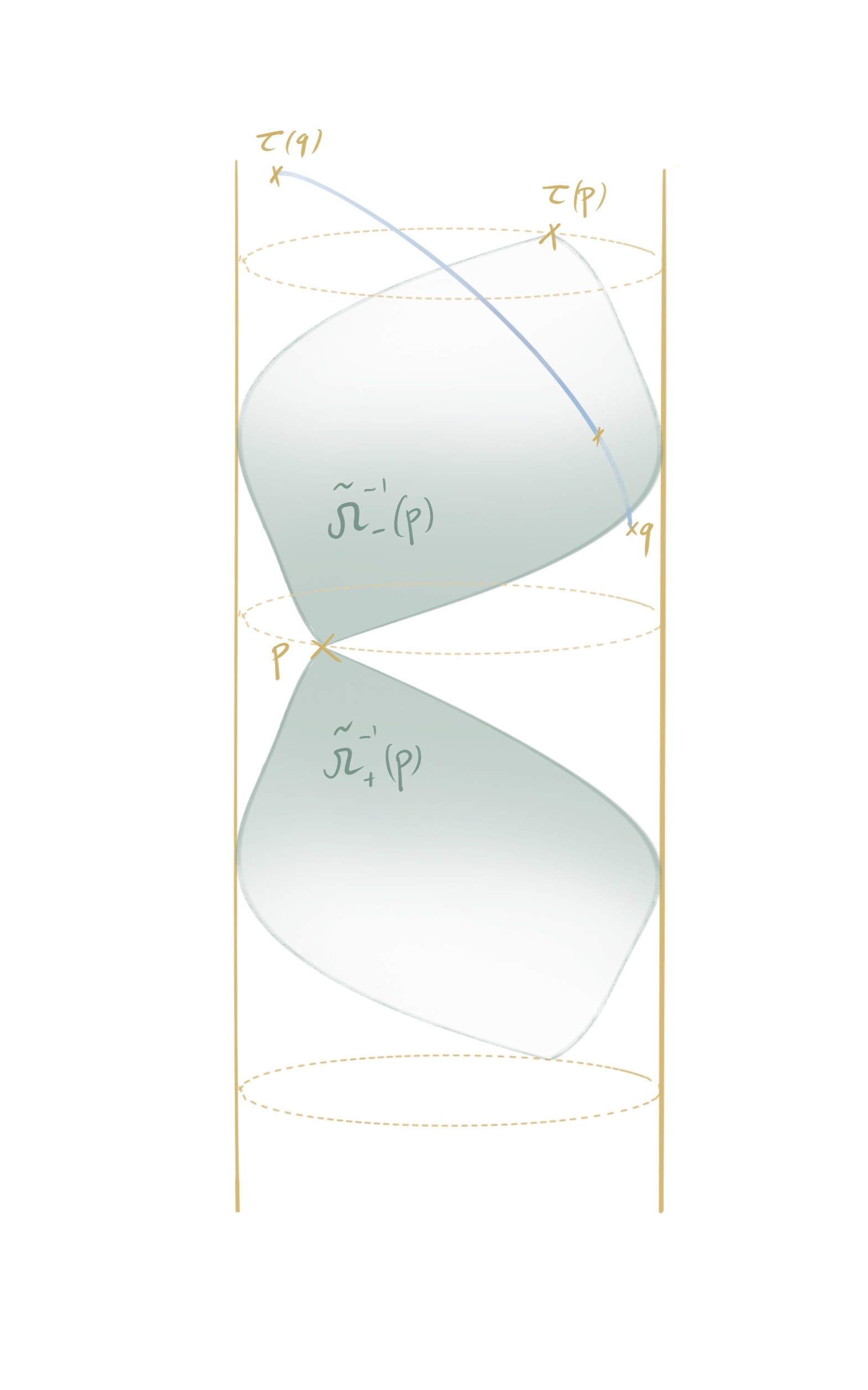} \centering
\caption{Any two complete totally geodesic lightlike hypersurfaces $\widetilde{\Omega}_-^{-1}(p), \widetilde{\Omega}_-^{-1}(q)\subseteq\widetilde{\AdS^{1,n}}$, emanating from the points $p,q\in \widetilde{\Ein}^n$ respectively, intersect if and only if $p$ and $q$ are not causally related.} \label{C_0_null_foliation} 
\end{figure}

\color{black}

\subsection{Codimension-one lightlike geodesic foliations of \texorpdfstring{$\widetilde{\AdS^{1,n}}$}{AdSn}}

In this section, we will describe all codimension-one lightlike geodesic foliations (or laminations) of $\widetilde{\AdS^{1,n}}$. Recall that any collection $X$ of complete totally geodesic light-like hypersurfaces is given by $\widetilde{\Omega}_-(X)$, a subset of $
\widetilde{\Ein}^{n}\simeq\mathbb R\times \S^{n-1}.$ By the same argument, a geodesic lightlike codimension-one lamination $\mathcal F$ of $\widetilde{\AdS^{1,n}}$ is a collection of pairwise disjoint connected complete geodesic lightlike hypersurfaces whose union is a closed subset.

\begin{example}\label{causalcurveconstruction}
Given an inextensible causal curve $\mathcal{C}$ in $\widetilde{\Ein}^n$, $\mathcal{C}$ is of the form $t\in \R \mapsto (t,c(t))\in \R\times\S^{n-1}$ with $c$ a unit speed geodesic in $\S^{n-1}$. Consider now for each $t\in\R$, the lightlike hypersurface $\widetilde{\Omega}_-^{-1}(\mathcal{C}(t))$. We claim that $(\widetilde{\Omega}_-^{-1}(\mathcal{C}(t)))_{t\in\R}$ defines a foliation of $\widetilde{\AdS^{1,n}}$. Notice that by Lemma \ref{intersection} and the fact that $\mathcal{C}$ is causal, the hypersurfaces do not intersect in $\widetilde{\AdS^{1,n}}$. In addition, if $(s,x)\in \widetilde{\AdS^{1,n}}$, consider the continuous function $t\in \R \mapsto t+d(c(t),x)- s\in \R$, it converges to $\pm \infty$ when $t$ does. Hence, there exists $t_0$ such that $t_0+d(c(t_0),x)=s$ which implies that $(s,x) \in \widetilde{\Omega}_-^{-1}(\mathcal{C}(t_0))$. We therefore conclude that $(\widetilde{\Omega}_-^{-1}(\mathcal{C}(t)))_{t\in\R}$ defines a codimension $1$ geodesic lightlike foliation of $\widetilde{\AdS^{1,n}}$.

\end{example}

\begin{exdef}[Rigid foliations] \label{exrigidfoliation}
Let $l\subset \R^{2,n}$ be an isotropic line. The hyperplane $l^\perp$ intersects $\AdS^{1,n}$ in a totally geodesic lightlike hypersurface, denoted $\F_l$. The hypersurface $\F_l$ has two connected components $\F_l^+$ and $\F_l^-$.

Two different hypersurfaces $\F_l$ and $\F_{l^\prime}$ intersect non-trivially if and only if the plane $l\oplus l^\prime$ has Lorentzian signature. 

Let now $P\subset
\R^{2,n}$ be an isotropic plane. The collection of hypersurfaces $\F_l$, $l\in \mathbb{P}(P)$, defines a foliation by lightlike geodesics hypersurfaces. We call a lightlike geodesic codimension-one foliation $\F$ of $\widetilde{\AdS^{1,n}}$ \textbf{rigid} if it is the inverse image of a foliation defined this way, i.e. corresponds to an isotropic plane $P\subset
\R^{2,n}$.

\end{exdef}
Notice that $\F$ is invariant under $-\Id$ and hence defines a codimension-one geodesic foliation in $\P\AdS^{1,n}$. In fact, these foliations are unique:
\begin{corollary}\label{rigidprojective}
    Any codimension-one geodesic foliation of $\P\AdS^{1,n}$ lifts to a rigid foliation of $\AdS^{1,n}$.
\end{corollary}
\begin{proof}
    Given a codimension-one geodesic lightlike foliation $\F$ of $\P\AdS^{1,n}$. It lifts to a foliation $\widehat{\F}$ of $\AdS^{1,n}$ invariant under $-\Id$. Now recall that any hyperplane $H\subset\R^{2,n}$ of signature $(0,-,+,..,+)$ intersects $\AdS^{1,n}$ in two disjoint codimension $1$ lightlike hypersurfaces, moreover, $-id$ interchange these hypersurfaces. This implies in particular that if one of them is a leaf of $\widehat{F}$, the second is a leaf too. Now, given two distinct hyperplanes $P$ and $Q$ of signatures $(0,-,+,..,+)$, they are defined as the orthogonals of two independent isotropic lines $\ell_P$ and $\ell_Q$. Consider the plane $L$ generated by $\ell_P$ and $\ell_Q$, suppose it is not isotropic i.e. $\ell_P$ and $\ell_Q$ are not orthogonal, now $P\cap Q$ is a codimension $2$ hyperplane orthogonal to both $\ell_P$ and $\ell_Q$, hence $P\cap Q \perp L$, in fact, since $L$ has signature $(-,+)$, we have $P\cap Q=L^\perp$ and in particular, it has signature $(-,+,...,+)$. This implies that $P\cap Q$ intersects $\AdS^{1,n}$ which is impossible, since it implies that the leaves defined by $P$ and $Q$ intersects. Thus $L$ is an isotropic plane and $\widehat{\F}$ is determined by $L$.
\end{proof}

    \color{black}
\begin{remark}[Space of rigid foliations]
    Recall that the space of (future-directed) lightlike geodesics of the globally hyperbolic spacetime $\R\times \S^{n-1}$ is a contact manifold which is contactomorphic to the unit co-sphere bundle (see e.g. \cite{Low2001, ChernovNemirovski2010}) $ST^*\S^{n-1}\simeq T^1\S^{n-1}$. Thus, $T^1\S^{n-1}$ is identified with the space of all rigid codimension-one geodesic lightlike foliations of $\widetilde{\AdS^{1,n}}$ on which the conformal group $\PO(2,n)$ acts by contactomorphisms.
\end{remark}

\begin{example}
    The unit bundle $T^1\S^1$ has two connected components corresponding to left and right-invariant codimension-one lightlike geodesic foliations of $\P\AdS^{1,2}=\PSL(2,\R)$.
\end{example}

\begin{theorem}\label{codimLight}
A collection $X$ of complete geodesic lightlike hypersurfaces defines a foliation if and only if $\widetilde{\Omega}_-(X)$ is an inextendible causal curve.
\end{theorem} 

\begin{proof}
The if part is given in example \ref{causalcurveconstruction}. Now suppose $X$ is a codimension-one lightlike geodesic foliation of $\widetilde{\AdS^{1,n}} \cong \R\times H_n$, this defines a canonical $1D$ null geodesic foliation $\mathcal{V}$ of $\widetilde{\AdS^{1,n}}$. Since $\gamma$ is timelike, it is transverse to $X$, and hence $\widetilde{\Omega}_-(X)$ is given by $\widetilde{\Omega}_-(\mathcal{V}(\gamma(t))_{t\in\R}$. By the intersection lemma \ref{intersection} and the foliation condition, this gives an inextendible causal curve in $ \widetilde{\Ein}^n$.

\end{proof}

 Unlike the projective anti-de Sitter space, codimension-one geodesic foliations of  $\AdS^{1,n}$ are not rigid. Indeed, by the theorem above, a geodesic foliation of $\AdS^{1,n}$ is given by $<\delta>$-invariant causal curve in $\widetilde{\Ein}^n$, thus, the projection of this causal curve in $\widehat{\Ein}^n$ is a closed degree $1$ loop. The rigidity in the projective model comes from the fact that a $\tau$-invariant causal curve is necessarily a lightlike geodesic. However, invariance by $\delta$ is not that restrictive.

In the $1+2$ dimension case, combining Theorem \ref{Dim1Codim1} and the theorem above, we get

\begin{corollary}
One-dimensional lightlike geodesic foliations of $\AdS^{1,2}$ are characterized by causal loops of degree $1$ in  $(\S^1 \times \S^1, dx^2-dy^2)$
\end{corollary}

\begin{remark}
    We note that any codimension-one geodesic foliation is locally Lipschitz (see \cite[Thm. 4]{zeghibGeodesicFoliationsofLorentz} for more details).
\end{remark}

\subsection{One-dimensional geodesic null foliations}

Observe that a codimension-one lightlike geodesic foliation yields a one-dimensional geodesic foliation. Indeed, each geodesic hypersurface is naturally foliated by null geodesics. In other words, we have the (saturation) map
\[
S:\{\text{co-dim 1 lightlike geod. foliations}\}
\longrightarrow
\{\text{1-d lightlike geod. foliations}\}.
\]
For $C^1$ regularity, Theorem \ref{Dim1Codim1} shows that, for $\widetilde{\AdS^{1,2}}$, the map $S$ is in fact bijective. In what follows, we show that the same conclusion holds for $C^0$ regularity.

\begin{theorem}\label{C0foliationofAdS}
    Let $\mathcal V$ be a one-dimensional geodesic lightlike $C^0$-foliation of $\widetilde{\AdS^{1,2}}$. Then the orthogonal distribution $T\mathcal V^{\perp}$ is integrable with totally geodesic leaves.
\end{theorem}

\begin{proof}

Let us show that the orthogonal leaf $l^\perp$ of a any geodesic $l\in\mathcal V$ is foliated by the geodesics of $\mathcal V$. 
Note that for any point $y$ in the boundary that lies in the causal domain determined by $x$, then any lightlike emanating from $y$ meet no lightlike geodesic that emanate from $x$ (Lemma \ref{intersection}).\ 

Hence, for any orthogonal leaf $l^\perp$ and for any lightlike geodesic $\tilde{l}$ that is transverse to it, their respective past limit points cannot be related by a causal curve. \

Let us consider $x\in\widetilde{\Omega}_-(\mathcal V) $ and $l_x\in\mathcal V$ a geodesic of the foliation $\mathcal V$ that emanates from $x$. We choose $y,z\in \widetilde{\Omega}_-(\mathcal V) $ with $l_y,l_z\in\mathcal V$ such that $l_y$ (respectively $l_z$) lies below (respectively above) $l_x^\perp$. If $l_x^\perp$ is not foliated by geodesics of $\mathcal V$ then there exists a curve $c:[0,1]\mapsto \widetilde{\AdS^{1,2}}$ that connects a point of $l_y$ to a point of $l_z$ and avoids the geodesics of $\mathcal V$. The curve $c$ descends to a curve $c_{\mathcal V}\subseteq \mathcal V$ that connects $l_y$ to $l_y$ and avoids all the geodesics belonging to $l_x^\perp$. By continuity, there exists a time $t\in[0,1]$ such that $c_{\mathcal V}(t) $ intersects $l_x$ which is a contradiction and thus completes the proof of the theorem.

\end{proof}

\color{black}

\begin{corollary}
    Geodesic lightlike one-dimensional $C^0$-foliations of the projective model $\PSL(2,\R)$ are rigid (i.e. left or right invariant).
\end{corollary} 
\begin{proof}
    By Theorem \ref{C0foliationofAdS}, any geodesic lightlike one-dimensional $C^0$-foliation $\mathcal{V}$ of $\AdS^{1,2}\cong \PSL(2,\R)$ gives rise to a geodesic lightlike codimension-one $C^0$-foliation $\F$ which is rigid by Corollary \ref{rigidprojective}. 
\end{proof}

\section{Timelike and mixed geodesic foliations of \texorpdfstring{$\AdS^{1,2}$}{H12}}
In this section, we study mixed geodesic foliations of the $ 3$-dimensional anti-de Sitter space $\AdS^{1,2}$. The projective model $\P\AdS^{1,2}$ can be identified (up to a scalar) with $\PSL(2,\R)$ endowed with the Killing form. Its double cover $\AdS^{1,2}$ is therefore identified with $\SL(2,\R)$. We get two descriptions of geodesics, either as cosets of $1$-parameter subgroups if we use the $(\P)\SL(2,\R)$ model or as intersections with $2$-planes if we consider the quadric model $\AdS^{1,2}$. 
\subsection{The Grassmannian identification}

Let $\Lambda^2 \R^{2,2}$ be the 6-dimensional space of $2$-forms. The scalar product $\langle \cdot, \cdot \rangle$, associated to $q$, on $\R^{2,2}$ induces a scalar product on $\Lambda^2 \R^{2,2}$ defined by:
$$ \langle v \wedge w, v' \wedge w' \rangle = \det \begin{pmatrix} \langle v, v' \rangle & \langle v, w' \rangle \\ \langle w, v' \rangle & \langle w, w' \rangle \end{pmatrix} $$
This scalar product is of signature $(4,2)$.

Given a plane $P$ of signature $(-,-)$ (respectively $(-,+)$) with orthonormal basis $\{v, w\}$, we can associate to it a unit $2$-form $\omega_P = v \wedge w$ which has norm $\langle \omega_P, \omega_P \rangle = 1$ ( respectively $-1$).

\subsubsection{The Hodge star splitting}
The Hodge star operator $* : \Lambda^2 \R^{2,2} \to \Lambda^2 \R^{2,2}$ can be defined using the scalar product above: $$ 
\text{for any} \ \alpha, \beta \in \Lambda^2 \R^{2,2}, \;\;\; \alpha \wedge * \beta = \langle \alpha, \beta \rangle \omega$$ 
 where $\omega$ is a fixed volume form (compatible with the orientation). The Hodge star satisfies $** = \Id$ in this case, this induces a decomposition of $\Lambda^2 \R^{2,2}$ into eigenspaces $E_\pm = \{ \omega : *\omega = \pm \omega \}$, where each $E_\pm$ inherits a metric of signature $(1,2)$ (the trace of the Hodge operator is $0$).
Any simple unit $2$-form $\omega_P$ decomposes uniquely as:
$$ \omega_P = \omega_+ + \omega_-, \quad \text{with } \omega_\pm \in E_\pm $$
The condition that $\omega_P$ is simple ($\omega_P \wedge \omega_P = 0$) and $\langle \omega_P, \omega_P \rangle = \pm 1$ implies:
$$ \langle \omega_+, \omega_+ \rangle =\pm \frac{1}{2} \quad \text{and} \quad \langle \omega_-, \omega_- \rangle = \pm \frac{1}{2} $$

Since the unit (+1) hyperboloid in a $(2,1)$ space can be identified with the hyperbolic plane $\H^2$, we have:
$$ G_2^{-}(\R^{2,2}) \cong \H^2 \times \H^2 $$
where $G_2^{-}(\R^{2,2})$ is the grassmannian of negative definite planes i.e. the space of timelike geodesics of $\AdS^{1,2}$. Similarly, for $(-,+)$ planes, we get: 
$$ G_2^{+}(\R^{2,2}) \cong \dS^{1,1} \times \dS^{1,1}$$
where $G_2^{+}(\R^{2,2})$ is the grassmannian of Lorentzian planes i.e. the space of spacelike geodesics of $\AdS^{1,2}$.
\subsubsection{The intersection lemma}
In order to study foliations, it would be crucial to determine when two geodesics intersect in $\AdS^{1,2}$. In the timelike case, this we get a nice criterion:

Given two oriented planes $P, Q$ in $\R^{2,2}$, seen as $2$-forms $\omega_P$ and $\omega_Q$ respectively. We can see directly that the intersection of $P$ and $Q$ is non trivial if and only if $\omega_P\wedge \omega_Q = 0$ (i.e. $\langle \omega_P, *\omega_Q \rangle = 0$). With this in mind, one can conclude
\begin{lemma}[Intersection condition]
Two oriented timelike planes $P(x, y)$ and $Q(x', y')$ in $\AdS^{1,2}$ (where $x, x', y, y' \in \H^2$) intersect non-trivially if and only if:
$$ d_{\H^2}(x, x') = d_{\H^2}(y, y') $$
They are disjoint if and only if $d_{\H^2}(x, x') \neq d_{\H^2}(y, y')$.
\end{lemma}

\begin{proof}
Two planes $P$ and $Q$ intersect non-trivially if and only if $\langle \omega_P, *\omega_Q \rangle = 0$. Using the decomposition $\omega_P = \omega_+ + \omega_-$ and $\omega_Q = \eta_+ + \eta_-$:
$$ \langle \omega_+ + \omega_-, *(\eta_+ + \eta_-) \rangle = \langle \omega_+ + \omega_-, \eta_+ - \eta_- \rangle = \langle \omega_+, \eta_+ \rangle - \langle \omega_-, \eta_- \rangle = 0 $$
In the hyperboloid model, $\langle a, b \rangle = \frac{1}{2}\cosh(d_{\H^2}(a, b))$. Thus:
$$ \frac{1}{2}\cosh(d_{\H^2}(x, x')) - \frac{1}{2}\cosh(d_{\H^2}(y, y')) = 0 \implies d_{\H^2}(x, x') = d_{\H^2}(y, y') $$
\end{proof}

\subsection{Foliations with timelike leaves}

Given a $C^0$ foliation $\mathcal{F}$ of $\AdS^{1,2}$ by complete geodesics. Assume moreover that $\mathcal{F}$ contains at least one timelike leaf. The set of points whose leaf is timelike is an open saturated set. Let $U\subset \AdS^{1,2}$ be a connected component of this set. Recall that the space of timelike geodesics in $\AdS^{1,2}$ can be identified with $\H^2 \times \H^2$. Therefore, we can define a map $\tau: U \rightarrow \H^2 \times \H^2 $ that associate to each point the leaf of $\F$ passing through it. We can derive the following characterization: 

\begin{proposition}
With the notations above and up to switching factors, $\tau(U)$ is the graph of distance-decreasing map
$f:\H^2\longrightarrow \H^2$.
Moreover, if $f$ is contraction then $\F$ is a global timelike foliation i.e. $U=\AdS^{1,2}$.
\end{proposition}

\begin{proof}
Since $\tau(U)$ defines a foliation, its leaves are pairwise disjoint. By the intersection lemma, for any $(x, y) \neq (x', y') \in \tau(U)$, $d(x, x') \neq d(y, y')$. Since $\tau(U)$ is a connected topological surface, the space $(\tau(U)\times\tau(U))\setminus \text{Diag}$ is connected. Hence, since the function $((x,y),(x',y'))\in (\tau(U)\times\tau(U))\setminus \text{Diag} \mapsto d(x,x')-d(y,y')$ 
is continuous and never vanishes, it cannot change sign. Therefore, we can suppose that $d(x,x')>d(y,y')$ for all distinct pairs $(x,y), (x',y')$ in $\tau(U)$. This
implies immediately that the projection on the $1$st factor $pr_1$ is injective. Since $\tau(U)$ is a surface, by invariance of domain the projection
$$
pr_1:\tau(U)\longrightarrow V:=pr_1(\tau(U))
$$
is a homeomorphism and $V$ is open in $\H^2$. Thus
$$
\tau(U)=\{(x,f(x)):x\in V\}
$$
for a continuous map $f:V\to\H^2$, moreover $d(f(x),f(x'))<d(x,x')$, for every $x\neq x'$.

Consider now $(x_n)$ a sequence in $V$ converging to some $x\in\H^2$. Since
$$
d(f(x_n),f(x_m))<d(x_n,x_m),
$$
the sequence $(f(x_n))$ is Cauchy. Since $\H^2$ is complete, there exists $y\in\H^2$ such that
$f(x_n)$ converges to $y$. Therefore
$(x_n,f(x_n))$ converges to $ (x,y)\in \H^2\times \H^2$. In particular, the leaves $\mathcal{C}_n$ corresponding to $(x_n,f(x_n))$ converge to the timelike geodesic $\mathcal{C}$ corresponding to $(x,y)$.

We claim that $\mathcal{C}$ is a leaf of $\F$. Choose $p\in\mathcal{C}$ and $p_n\in\mathcal{C}_n$ such that $p_n$ converges to $p$. Consider a trivialization neighborhood for foliation around $p$. Since $\F$ is a $C^0$ foliation, the leaves through $p_n$ converge to the leaf through $p$. On the other hand, since $\mathcal{C}_n$ converges to $\mathcal{C}$, these leaves converge to a piece of $\mathcal{C}$. Hence the leaf through $p$ coincides locally with $\mathcal{C}$ and, since both are geodesics, it is equal to $\mathcal{C}$. Thus $\mathcal{C}\in\F$.

Moreover, $\mathcal{C}$ is timelike. Since $p_n\in U$, $p_n$ converges to $p$, and $U$ is a connected component of the set of points whose leaf is timelike, we have $p\in U$. Since $U$ is saturated, $\mathcal{C}\subset U$. Therefore $(x,y)\in\tau(U)$ and $x\in V$. We conclude that $V$ is closed in $\H^2$. Since $V$ is non-empty, open and closed in $\H^2$, we have $V=\H^2$. 

Suppose Now that the map $f$ is contracting i.e. there exists $\lambda<1$, such that 
$$d(f(x),f(y)) \leq \lambda d(x,y)$$
for any $x,y\in \H^2$. Then, the graph of $f$ intersects the graph of any element in $\PSL(2,\R)$. Indeed, for any $g\in\PSL(2,\R)$, the map $g^{-1}\circ f$ admits a fixed point $x\in\H^2$ (by Banach fixed-point theorem). Hence, the point $(x,f(x))$ is in the graph of $g$. In particular, $\F$ admits a leaf that passes through $g\in \PSL(2,\R)\cong \AdS^{1,2}$. Since this leaf belongs to $\tau(U)$, we have $g\in U$. This implies that $\F$ is a global timelike foliation.
\end{proof}

\begin{remark}
    Note that this characterization remains the same for all coverings of $\AdS^{1,2}$, since covering automorphisms preserve individually each timelike geodesic.
\end{remark}

\subsection{Some examples}

The first simple example of a global timelike foliation can be given by taking right cosets of an elliptic $1$-parameter subgroup in  $\PSL(2,\R)$, This foliation is given by a constant map on $\H^2$. It is worth nothing, that not every distance-decreasing map $f: \H^2 \rightarrow \H^2$ defines a global timelike foliation. Indeed, we can construct mixed foliations of $\AdS^{1,2}$ by restricting global foliations of $\S^3$ on the image of $\AdS^{1,2} \subset \S^3$. 
\subsubsection{\textbf{Mixed foliations (TLS)}}
Recall that  $\AdS^{1,2}$ is identified with the solid torus $ \mathcal A=\{x\in \S^3\mid q(x)<0\}$ in $\S^3 \subset \C^2$. Moreover, the boundary $\widehat{\Ein}^2$ is identified with the Clifford torus $ \mathcal T=\S^1\left(\frac1{\sqrt2}\right)\times \S^1\left(\frac1{\sqrt2}\right).$ Geodesics in this case, are projectivisation of $2$-planes. Hence, all anti-de Sitter geodesics lie on great circles of $\S^3$. Their causal type can be read directly from the way the great circle intersects the Clifford torus $\mathcal T$.

\begin{proposition}
Let $C_P=P\cap \S^3$ be the great circle associated to a geodesic.

\begin{enumerate}
\item If $C_P\cap\mathcal T=\varnothing$, then $q_{\mid P}$ is negative definite and $C_P\subset\mathcal A$. The corresponding geodesic is timelike.

\item If $C_P$ is tangent to $\mathcal T$, then $q_{\mid P}$ is degenerate. The circle is tangent to $\mathcal T$ at two antipodal points, and the connected components of $C_P\setminus\mathcal T$ contained in $\mathcal A$ are lightlike geodesics.

\item If $C_P$ intersects $\mathcal T$ transversely, then $q_{\mid P}$ has signature $(-,+)$. The circle intersects $\mathcal T$ in four points and alternates between the regions $q<0$ and $q>0$. The two connected components of $C_P\cap\mathcal A$ are spacelike geodesics.
\end{enumerate}
\end{proposition}

We are now ready to construct mixed foliations. Let $\mathcal{F}$ be a foliation of $\S^3$ by great circles, then its restriction to $\mathcal{A}$ gives a geodesic foliation of $\AdS^{1,2}$, where each leaf is a connected component of the intersection of a great circle of $\mathcal{F}$ with $\mathcal{A}$. The Hopf fibration gives a first example. Indeed, the Hopf fibers are great circles, and since they come from complex lines in $\C^2$, the function $q(z_1,z_2)=-|z_1|^2+|z_2|^2$ is constant along each Hopf fiber. Hence every Hopf fiber which intersects $\mathcal{A}$ is entirely contained in $\mathcal{A}$, and the Hopf fibration induces a global timelike geodesic foliation of $\AdS^{1,2}$.

More generally, applying an element of $\PGL(4,\R)$ to the Hopf fibration gives another foliation of $\S^3$ by great circles. For a general projective transformation, the resulting great circles may intersect the Clifford torus $\mathcal{T}$ transversely, and the induced foliation of $\mathcal{A}$ may therefore contain leaves of different causal types: timelike, lightlike and spacelike. Notice that if the projective transformation is $\C$-linear, it preserves the family of complex lines and hence simply permutes the Hopf fibers; thus one needs a totally real projective deformation in order to obtain new causal behavior.

\subsubsection{\textbf{Timelike-Lightlike foliations}}

The previous construction produces either global timelike foliations or mixed foliations with leaves of each causal type, It can not however produce mixed foliations containing timelike and lightlike leaves but no spacelike leaves. Indeed, given a foliation $\F$ of $\S^3$ by great circles, consider $V$ a unit vector field generating it. The restriction of $\F$ on $\mathcal{A}$ has no spacelike leaves if and only if whenever a leaf of $\F$ intersects the Clifford torus $\mathcal{T}$, it is tangent to it. In particular, $V_{\mid \mathcal{T}}$ is tangent to $\mathcal{T}$, and this implies that any leaf of $\F$ tangent to $\mathcal{T}$ is necessarily contained in it. Thus, the restriction of $\F$ on $\mathcal{A}$ has only timelike leaves. 

Nevertheless, we can still construct geodesic foliations of $\AdS^{1,2}$ that have only timelike and lightlike leaves. For this, we will use the projective model. Recall that the quotient map $\S^3\to\mathbb{RP}^3$ identifies antipodal points and sends a great circle $C_P=P\cap \S^3$ to the projective line $\mathbb P(P)\subset\mathbb{RP}^3$. Thus, in the projective model, the geodesics of $\P\AdS^{1,2}$ are the connected components of
$$ \mathbb P(P)\cap\mathbb P(\{q<0\}). $$
Recall that an affine chart of $\mathbb{RP}^3$ is a complement of a projective plane. More precisely, if $\ell:\R^4\to\R$ is a non-zero linear form, then
$$ U_\ell=\{[x]\in\mathbb{RP}^3\mid \ell(x)\neq0\} $$
is identified with the affine hyperplane
$$ \{x\in\R^4\mid \ell(x)=1\}\simeq\R^3 $$
by sending $[x]$ to $x/\ell(x)$. Under this identification, every projective line which is not contained in the plane at infinity $\{\ell=0\}$ becomes an ordinary affine line. Hence, in a suitable affine chart, anti-de Sitter geodesics are represented by straight lines, or by half lines inside the domain $\mathbb P(\{q<0\})$. The causal type is then given by the intersection with the boundary quadric $\mathbb P(\{q=0\})$: a timelike geodesic corresponds to a projective line which does not meet the boundary, a lightlike geodesic to a projective line tangent to the boundary, and a spacelike geodesic to a projective line intersecting the boundary transversely in two points.

Consider an affine chart $U=\{[x_0:x_1:x_2:x_3]\in \R\mathbb{P}^3, x_0\neq 0 \}  $ identified with $\R^3$. As seen before, $\P\AdS^{1,2}\cap U$ is identified with the interior of the open "ball" $B_1=\{(x,y,z) \in \R^3, \; x^2+y^2-z^2 < 1\}$. Moreover, $\P\AdS^{1,2} \setminus B_1$ is identified with the open disc $D=\{[0:1:x_2:x_3]\in \RP^3 , x_2^2+x_3^2 <1\}$. Now, we consider a smooth function $f:]0,1[ \rightarrow \R$ such that $f(r)=1$ when $r\in[1/2,1[$ and is strictly decreasing on $]0,1/2]$ with limit $0$ at $0$. The ball $B_1=\bigsqcup_{0<r<1} S_r \sqcup \{x=y=0\}$, where $S_r=\{x^2+y^2-f(r)z^2=r\}$. To obtain our foliation, notice that the surfaces $S_r$ are ruled surfaces and can be naturally foliated by lines (smoothly with respect to $r$). In addition, we notice that these lines are asymptotic to the boundary of $\P\AdS^{1,2}\cap U$ when $r>1/2$ and hence are lightlike lines. On the other hand, when $0<r<1/2$ these lines are timelike and their intersection with $D$, covers bijectively $D\setminus \{0\}$. Finally, the line $\{x=y=0\}$ is a timelike geodesic of $\P\AdS^{1,2}$ that intersects $D$ in the origin. In conclusion, this construction provides a smooth foliation of $\P\AdS^{1,2}$ by timelike and lightlike geodesics only. In addition, this foliation is invariant by the action $\SO(2)$ on $\P\AdS^{1,2}$ by conjugacy. (see Figure \ref{T-L_foliations} below).

\begin{figure}[h] 
\includegraphics[scale=0.6]{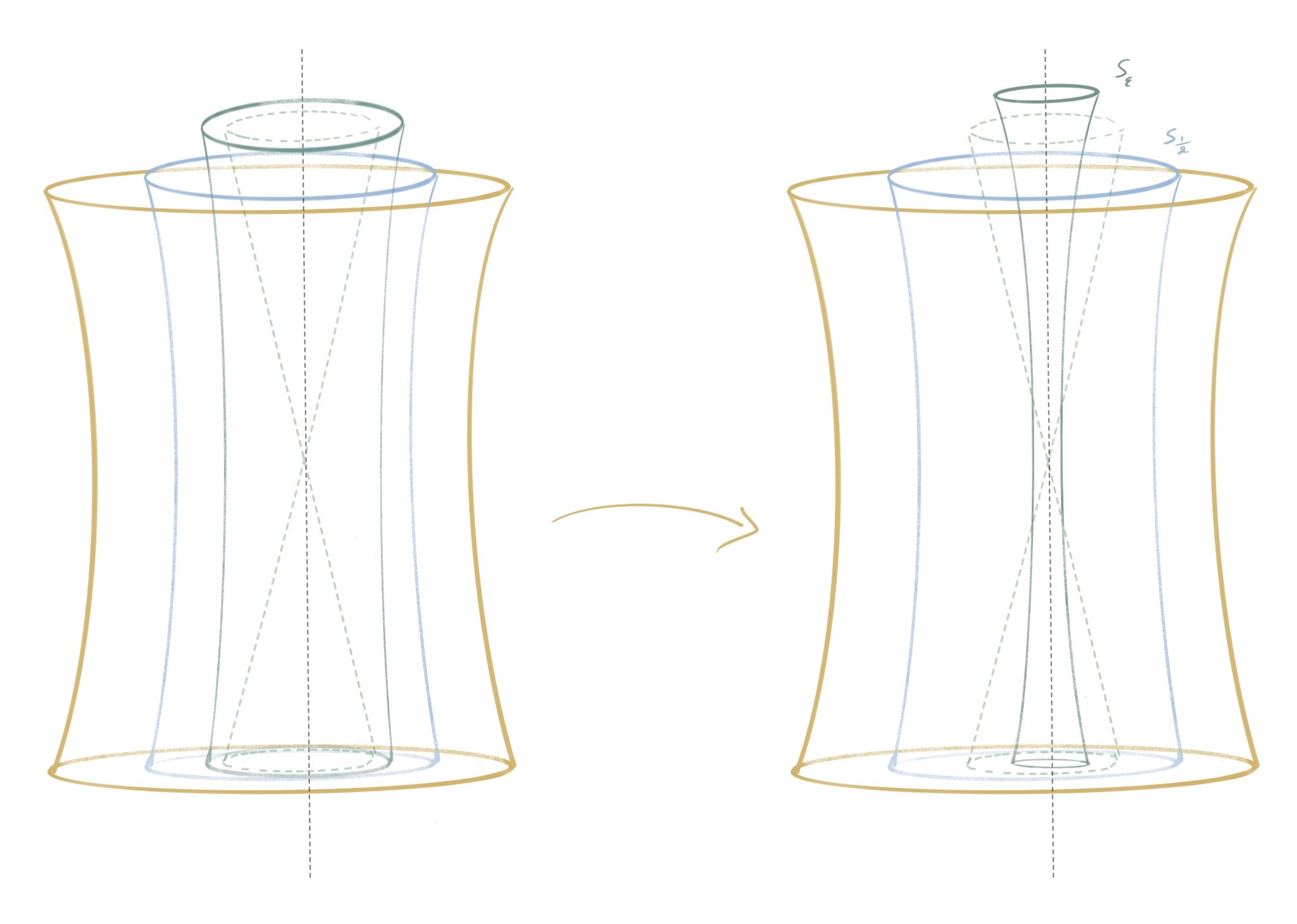} \centering
\caption{The family of hyperboloid $H_r$ of equations $\{x^2+y^2-z^2=r\}$ where $r\in]0,1[$, foliates the complementary region of the cone $C$, of equation $\{x^2+y^2-z^2=0\}$, in $\P\AdS^{1,2}\cap U$. The aforementioned construction is given by deforming continuously and bijectively the region lying between the hyperboloid $H_{\frac{1}{2}}$ and the cone $C$ onto the region between $H_{\frac{1}{2}}$ and the $z$-axis, by dilating each hyperboloid $H_r$ along the latter by a factor that depends on $r$.} \label{T-L_foliations} 
\end{figure}

\subsubsection{\textbf{Spacelike and Spacelike-Lightlike foliations}}
Consider an isotropic line $\ell$ in $\R^{2,n}$ and consider the family of planes $(P_i)_{i\in I}$ containing $\ell$ which are either of signature $(0,-)$ or $(+,-)$. This family defines a geodesic foliation of $\AdS^{1,n}$ with spacelike and lightlike leaves. Indeed, for any point $p\in \AdS^{1,n}$, the plane generated by $p$ and $\ell$ is in the family, on the other hand, for any $j, i \in I$, the intersection $P_i \cap P_j=\ell$ if and only if $i\neq j$. Hence, the intersections $\F_i=P_i\cap \AdS^{1,n}$ is a geodesic foliation. Moreover, it is analytic. Notice that this foliation is characterized by a unique point $x$ in the boundary and is the only $ \mathcal {P}_1$-invariant foliation, where $\mathcal{P}_1$ is the stabilizer of $x$ in $\PO(2,n)$.

One way to obtain geodesic spacelike foliations of $\P\AdS^{1,2}\cong\PSL(2,\R)$ is by taking the right cosets of a $1$-parameter spacelike subgroup $H$ of $\PSL(2,\R)$. These foliations are invariant under the entire right action of $\PSL(2,\R)$. Moreover, any geodesic foliation of the hyperbolic plane $\H^2$ gives rise to an $\P\SO(2)$-invariant spacelike foliation of $\P\AdS^{1,2}$ (which can be lifted to $\AdS^{1,2}$). Indeed, Given a geodesic foliation $\F$ of $\H^2$, since the projection $\PSL(2,\R)\rightarrow \H^2\cong\P\SO(2)\backslash\PSL(2,\R)$ is a pseudo-Riemannian submersion with fibers $\P\SO(2)g$, any leaf of $\F$ lifts to a unique spacelike geodesic of $\AdS^{1,2}$ once a point in the fiber is fixed. We can verify that all the lifts of all the leaves of $\F$ define a foliation (by the uniqueness of the horizontal lifting). In particular, the space of spacelike geodesic foliations contains a copy of the space of geodesic foliations of $\H^2$ which can be identified with causal inextendible curves in $\dS^{1,1}$ (see §3 in \cite{zeghibLipschitzregularity}). This flexibility of spacelike geodesic foliation is no longer true for compact $\AdS^{1,2}$-manifolds \cite{geodfoliationsofcompact}.

\subsection{Local timelike foliations of even dimensional anti-de Sitter spaces}

A classical result \cite[Lemma 2]{Ferusfoliationsofspheres} proves that even-dimensional spheres do not admit local foliation by great circles. This implies in particular that even-dimensional anti-de Sitter spaces do not admit local timelike foliations (this has already been proved in  {\cite[Cor.8.3]{Tholozan_timelike_foliation}} by another method). Indeed, the quadric model $\AdS^{1,n}$ can be seen as an open set of $\S^n$, using projection by positive scalars. Moreover, timelike geodesics are great circles. Thus, any local foliation of $\AdS^{1,2n-1}$ by complete timelike geodesics is a local foliation of $\S^{2n}$ by great circles. In fact, we can apply the same idea on any anti-de Sitter space $\AdS^{p,q}$ of signature $(p,q)$, i.e., the $-1$ ball in $\R^{p+1,q}$, provided $p+q$ is even and $p\geq 2$. Indeed, in this case $\AdS^{p,q}$ can be identified with an open subset $\mathcal{A}^{p,q}$ of the unit sphere $\S^{p+q}$ and timelike geodesics of $\AdS^{p,q}$ are great circles lying entirely in $\mathcal{A}^{p,q}$.

For the sake of completeness, we give here a proof of this result:

\begin{proposition}[\cite{Ferusfoliationsofspheres}]
    The sphere $\S^{2n}$ does not admit a local foliation by great circles. 
\end{proposition}

\begin{proof}
    Let $\F$ be a foliation of an open set $U\subset \S^{2n}$. Consider $X\in\mathcal{X}(U)$, a unit vector field generating $\F$ and $\varphi_t$ its flow. 

    Since $\F$ is a geodesic foliation, we get that for any $p\in U$ and any $v\in T_pU$, the vector field $J(t):=d_p\varphi_t(v)$ defined along the orbit $\varphi_t(p))_{t\in \R}$ is a Jacobi field. In particular, since $U\subset \S^{2n}$ has constant curvature $1$, we have $$J(t)=P_t(\cos(t)v+\sin(t) A(v))$$
    where $P_t$ is the parallel transport along $\varphi_t(p))_{t\in \R}$ and $A: w\in T_pU \mapsto \nabla_w X \in T_pU$ for the Levi-Civita connection $\nabla$. Since $X$ is has constant norm $1$ and $\nabla$ is compatible with the round metric on $U$, we get that $A$ preserves the orthogonal complement $X^\perp$ of $X$ at each point. 

    Now, since the restriction of $A$ on $X^\perp$ is a linear endomorphism of a linear space of dimension $2n-1$. Hence, there exists a unit vector $v\in X^\perp$ and $\lambda in \R$, such that $A(v)=\lambda v$. In particular, the associated Jacobi field has the following form
    $$ J(t)=(\cos(t)+\lambda \sin(t))P_t(v)$$
    for all $t\in \R$. Notice that, we can find $t_0\in\R$ for which $\cos(t_0)+\lambda \sin(t_0)=0$, therefore $J(t_0)=0$ which is impossible since $J(t_0)=d\varphi_{t_0}(v)\neq 0$ since $\varphi_{t_0}$ is a diffeomorphism of $U$. We conclude that the foliation $\F$ cannot exist.
    
\end{proof}

\bibliographystyle{plain}
\bibliography{Bib.bib}

\end{document}